\documentclass[10pt,twoside]{amsart}

\usepackage{amsthm}
\usepackage{amsmath}
\usepackage{amssymb}
\usepackage{amsfonts}
\usepackage{latexsym}
\usepackage{enumitem}
\usepackage{mathrsfs}
\usepackage[all]{xy} \SelectTips{eu}{}
\usepackage{hyperref}

\newtheorem*{intsetup*}{Setup}

\newtheorem*{intNotations*}{Notations and Conventions}

\newtheorem{intthm}{Theorem}[]

\newtheorem*{intque*}{Question}
\newtheorem*{intexa*}{Example}

\newcommand{\numberseries}{\bfseries}   
\newlength{\thmtopspace}                
\newlength{\thmbotspace}                
\newlength{\thmheadspace}               
\newlength{\thmindent}                  
\newtheoremstyle{bfupright head,slanted body}
                {\thmtopspace}{\thmbotspace}
                {\slshape}{\thmindent}{\bfseries}{.}{\thmheadspace}
                {{\numberseries \thmnumber{#2\;}}\thmnote{#3}}

\newtheoremstyle{bfupright head,upright body}
                {\thmtopspace}{\thmbotspace}
                {\upshape}{\thmindent}{\bfseries}{.}{\thmheadspace}
                {{\numberseries \thmnumber{#2\;}}\thmnote{#3}}

\newtheoremstyle{fixed bf head,slanted body}
                {\thmtopspace}{\thmbotspace}{\slshape}
                {\thmindent}{\bfseries}{.}{\thmheadspace}
                {{\numberseries \thmnumber{#2\;}}\thmname{#1}\thmnote{ (#3)}}

\newtheoremstyle{fixed bf head,upright body}
                {\thmtopspace}{\thmbotspace}{\upshape}
                {\thmindent}{\bfseries}{.}{\thmheadspace}
                {{\numberseries \thmnumber{#2\;}}\thmname{#1}\thmnote{ (#3)}}

\newtheoremstyle{numbered paragraph}
                {\thmtopspace}{\thmbotspace}{\upshape}
                {\thmindent}{\upshape}{}{\thmheadspace}
                {{\numberseries \thmnumber{#2.}}}

\theoremstyle{bfupright head,slanted body}
\newtheorem{res}{}[section]             \newtheorem*{res*}{}

\theoremstyle{bfupright head,upright body}
\newtheorem{bfhpg}[res]{}               \newtheorem*{bfhpg*}{}

\theoremstyle{fixed bf head,slanted body}
\newtheorem{thm}[res]{Theorem}          \newtheorem*{thm*}{Theorem}
\newtheorem{prp}[res]{Proposition}      \newtheorem*{prp*}{Proposition}
\newtheorem{cor}[res]{Corollary}        \newtheorem*{cor*}{Corollary}
\newtheorem{lem}[res]{Lemma}            \newtheorem*{lem*}{Lemma}
         \newtheorem*{que*}{Question}
\theoremstyle{fixed bf head,upright body}
\newtheorem{setup}[res]{Setup}           \newtheorem*{setup*}{Setup}
\newtheorem{dfn}[res]{Definition}       \newtheorem*{dfn*}{Definition}
\newtheorem{rmk}[res]{Remark}           \newtheorem*{rmk*}{Remark}
\newtheorem{exa}[res]{Example}           \newtheorem*{exa*}{Example}
           \newtheorem*{exer*}{Exercise}

\theoremstyle{numbered paragraph}
\newtheorem{ipg}[res]{}

\newlength{\thmlistleft}        
\newlength{\thmlistright}       
\newlength{\thmlistpartopsep}   
\newlength{\thmlisttopsep}      
\newlength{\thmlistparsep}      
\newlength{\thmlistitemsep}     

\newcounter{eqc}
\newenvironment{eqc}{\begin{list}{\upshape (\textit{\roman{eqc}})}%
    {\usecounter{eqc}%
      \setlength{\leftmargin}{\thmlistleft}%
      \setlength{\labelwidth}{\thmlistleft}%
      \setlength{\rightmargin}{\thmlistright}%
      \setlength{\partopsep}{\thmlistpartopsep}%
      \setlength{\topsep}{\thmlisttopsep}%
      \setlength{\parsep}{\thmlistparsep}%
      \setlength{\itemsep}{\thmlistitemsep}}}%
  {\end{list}}%

\newcounter{prt}
\newenvironment{prt}{\begin{list}{\upshape (\alph{prt})}%
    {\usecounter{prt}%
      \setlength{\leftmargin}{\thmlistleft}%
      \setlength{\labelwidth}{\thmlistleft}%
      \setlength{\rightmargin}{\thmlistright}%
      \setlength{\partopsep}{\thmlistpartopsep}%
      \setlength{\topsep}{\thmlisttopsep}%
      \setlength{\parsep}{\thmlistparsep}%
      \setlength{\itemsep}{\thmlistitemsep}}}%
  {\end{list}}%

\newcounter{rqm}
\newenvironment{rqm}{\begin{list}{\upshape (\arabic{rqm})}%
    {\usecounter{rqm}%
      \setlength{\leftmargin}{\thmlistleft}%
      \setlength{\labelwidth}{\thmlistleft}%
      \setlength{\rightmargin}{\thmlistright}%
      \setlength{\partopsep}{\thmlistpartopsep}%
      \setlength{\topsep}{\thmlisttopsep}%
      \setlength{\parsep}{\thmlistparsep}%
      \setlength{\itemsep}{\thmlistitemsep}}}%
  {\end{list}}%

  {\end{list}}%

\newenvironment{prf*}[1][Proof]{%
  \begin{proof}[\bf #1]
    \setcounter{equation}{0}
    }
  {\end{proof}
}

\newcommand{\pgref}[1]{\ref{#1}}

\renewcommand{\eqref}[1]{(\pgref{eq:#1})}

\newcommand{\pgcite}[2][?]{\cite[Page~#1]{#2}}
\newcommand{\thmcite}[2][?]{\cite[Theorem~#1]{#2}}
\newcommand{\prpcite}[2][?]{\cite[Proposition~#1]{#2}}

\newcommand{\corcite}[2][?]{\cite[Corollary~#1]{#2}}
\newcommand{\lemcite}[2][?]{\cite[Lemma~#1]{#2}}

\newcommand{\dfncite}[2][?]{\cite[Definition~#1]{#2}}

\newcommand{\exacite}[2][?]{\cite[Example~#1]{#2}}

\newcommand{\eqclbl}[1]{{\upshape(\textit{#1})}}
\newcommand{\proofofimp}[3][:]{\mbox{\eqclbl{#2}$\!\implies\!$\eqclbl{#3}#1}}
\numberwithin{equation}{res}

\def\urltilda{\kern -.15em\lower .7ex\hbox{\~{}}\kern .04em}

\newcommand{\GF}{\mathsf{GF}}
\newcommand{\GP}{\mathsf{GP}}

\newcommand{\GI}{\mathsf{GI}}

\newcommand{\tr}{\mathsf{tr}}

\newcommand{\sfT}{\mathsf{T}}

\newcommand{\Proj}{\mathsf{Proj}}
\newcommand{\id}{\mathrm{id}}
\newcommand{\fd}{\mathrm{fd}}
\newcommand{\pd}{\mathrm{pd}}
\newcommand{\Inj}{\mathsf{Inj}}
\newcommand{\Flat}{\mathsf{Flat}}
\newcommand{\Mod}{\mathsf{Mod}}

\newcommand{\Coind}{\mathrm{Coind}}
\newcommand{\Ind}{\mathrm{Ind}}
\newcommand{\Id}{\mathrm{Id}}

\newcommand{\xra}[2][]{\xrightarrow[#1]{\:#2\:}}

\newcommand{\Rop}{R^{\sf op}}
\newcommand{\op}{\sf op}

\newcommand{\cok}{\mbox{\rm Coker}}
\newcommand{\kernel}{\mbox{\rm Ker}}
\newcommand{\im}{\mbox{\rm Im}}

\newcommand{\Hom}{\operatorname{Hom}}

\newcommand{\Ext}{\operatorname{Ext}}
\newcommand{\Tor}{\operatorname{Tor}}

\newcommand{\is}{\cong}

   \def\soft#1{\leavevmode\setbox0=\hbox{h}\dimen7=\ht0\advance
    \dimen7 by-1ex\relax\if t#1\relax\rlap{\raise.6\dimen7
    \hbox{\kern.3ex\char'47}}#1\relax\else\if T#1\relax
    \rlap{\raise.5\dimen7\hbox{\kern1.3ex\char'47}}#1\relax
    \else\if d#1\relax\rlap{\raise.5\dimen7\hbox{\kern.9ex
    \char'47}}#1\relax\else\if D#1\relax\rlap{\raise.5\dimen7
    \hbox{\kern1.4ex\char'47}}#1\relax\else\if l#1\relax
    \rlap{\raise.5\dimen7\hbox{\kern.4ex\char'47}}#1\relax
    \else\if L#1\relax\rlap{\raise.5\dimen7\hbox{\kern.7ex
    \char'47}}#1\relax\else\message{accent \string\soft
    \space #1 not defined!}#1\relax\fi\fi\fi\fi\fi\fi}

\begin{document}

\title[Gorenstein homological modules over tensor rings]%
{Gorenstein homological modules over tensor rings}

\author[Z.X. Di]{Zhenxing Di}
\address{Zhenxing Di: School of Mathematical Sciences,
Huaqiao University, Quanzhou 362000, China}
\email{dizhenxing@163.com}

\author[L. Liang]{Li Liang}
\address{Li Liang (Corresponding author):
Department of Mathematics, Gansu Center for Fundamental Research in Complex Systems Analysis and Control, Lanzhou Jiaotong University, Lanzhou 730070, China}
\email{lliangnju@gmail.com}
\urladdr{https://sites.google.com/site/lliangnju}

\author[Z.Q. Song]{Zhiqian Song}
\address{Zhiqian Song: School of Mathematical Sciences,
Huaqiao University, Quanzhou 362000, China}
\email{limbooosong@163.com}

\author[G.L. Tang]{Guoliang Tang}
\address{Guoliang Tang: School of Mathematical Sciences,
Zhejiang Normal University, Jinhua 321004, China}
\email{tangguoliang970125@163.com}

\thanks{Z.X. Di was partly supported by NSF of China (Grant No. 12471034) and
the Scientific Research Funds of Huaqiao University (Grant No. 605-50Y22050);
L. Liang was partly supported by NSF of China (Grant No. 12271230) and the Foundation for Innovative Fundamental Research Group Project of Gansu Province (Grant No. 25JRRA805).}


\keywords{Tensor ring; Gorenstein projective module; Gorenstein injective module;
Gorenstein flat module.}

\footnotetext{2020 \emph{Mathematics Subject Classification}. 18G25; 16D90.}


\begin{abstract}
For a tensor ring $T_R(M)$,
under certain conditions, we characterize the Gorenstein projective modules over $T_R(M)$, and prove that a $T_R(M)$-module $(X,u)$ is Gorenstein projective if and only if $u$ is monomorphic and $\cok(u)$ is a Gorenstein projective $R$-module.
Gorenstein injective (resp., flat) modules over $T_R(M)$ are also explicitly described.
Moreover, we give a characterization for the coherence of $T_R(M)$.
Some applications to trivial ring extensions and Morita context rings are given.
\end{abstract}

\maketitle

\thispagestyle{empty}

\section*{Introduction}
\label{Preliminaries}
\noindent
Throughout the paper, all rings are nonzero associative rings with identity and
all modules are unitary.
For a ring $R$, we adopt the convention that an $R$-module is a left $R$-module;
right $R$-modules are viewed as modules over the opposite ring $\Rop$.

The concept of the Gorenstein dimension of finitely generated modules over noetherian rings was introduced by Auslander and Bridger \cite{1969AB}.
It has been extended to modules over any associative rings by means of the notion of Gorenstein projective modules introduced by Enochs and Jenda \cite{1995EnochsGP}.
Moreover, they introduced the notion of Gorenstein injective modules in \cite{1995EnochsGP}, and the notion of Gorenstein flat modules together with Torrecillas in \cite{GF1993}.
One of the most important tasks in this topic is to describe
Gorenstein homological modules over a given extended ring.
In a series of papers \cite{2020GPTRI, TriExtGPMao, GFGN, 2022XI},
Gorenstein homological modules over triangular matrix rings,
trivial ring extensions and Morita context rings are investigated.

Let $R$ be a ring and $M$ an $R$-bimodule.
Recall from \cite{1991Cohn, BCALG} that
the tensor ring $T_R(M)$ is defined as $T_R(M)=\bigoplus_{i=0}^\infty M^{\otimes_Ri}$,
where $M^{\otimes_R0}=R$ and
$$M^{\otimes_R(i+1)}=M\otimes_R(M^{\otimes_Ri})$$
for $i \geqslant 0$.
We mention that the addition of elements in $T_R(M)$ is componentwise,
and the multiplication of elements in $T_R(M)$ is induced by the isomorphism
\begin{align*}
(\bigoplus_{n=0}^\infty M^{\otimes_Rn})\otimes_R(\bigoplus_{m=0}^\infty M^{\otimes_Rm})
\cong \bigoplus_{n,m=0}^\infty M^{\otimes_R(n+m)}.
\end{align*}
Examples of tensor rings include but are not limited to
trivial ring extensions, Morita context rings and triangular matrix rings.

Suppose that $M$ is $N$-\emph{nilpotent}, that is,
$M^{\otimes_R(N+1)}=0$ for some $N\geqslant 0$. Then the tensor ring $T_R(M)$ has an advantage that a module over $T_R(M)$ can be written as a pair $(X,u)$
with $X$ an $R$-module and $u \in \Hom_R(M\otimes_RX,X)$.
The classical homological properties of the tensor ring $T_R(M)$ are studied
in \cite{2012Ample, 1975Roganov}.
Recently, Chen and Lu \cite{TRchen} studied the Gorenstein
homological properties of modules over the tensor ring $T_R(M)$ for
an $N$-nilpotent $R$-bimodule $M$, where $R$ is a noetherian ring, and $M$ is finitely generated on both sides. They characterized finitely generated Gorenstein projective $T_R(M)$-modules.
To be more specific, they proved that under certain conditions,
a module over $T_R(M)$, that is, a pair $(X,u)$,
is finitely generated Gorenstein projective if and only if $u$ is monomorphic
and $\cok(u)$ is a finitely generated Gorenstein projective $R$-module.
For an arbitrary associative ring $R$, the following natural question arises:

\begin{que*} How to describe (infinitely generated) Gorenstein projective,
injective and flat modules over $T_R(M)$ for an $N$-nilpotent $R$-bimodule $M$?
\end{que*}

The main purpose of the present paper is to provide an answer for the above question and further give some relevant applications.
For a subcategory $\sf{X}$ of $R$-modules
and a subcategory $\sf{Y}$ of $R^{\op}$-modules,
the symbols $\Phi(\sf{X})$ and $\Psi(\sf{Y})$
denote the subcategories of $\Mod(T_R(M))$ and $\Mod(T_R(M)^{\op})$, respectively,
which are defined as follows:
\begin{align*}
&\Phi({\sf{X}}) = \{(X,u) \in \Mod(T_R(M)) \,|\,
u \textrm{ is a monomorphism and } \cok(u) \in \sf{X}\}; \\
&\Psi({\sf{Y}}) = \{[Y,v] \in \Mod(T_R(M)^{\op}) \,|\,
v \textrm{ is an epimorphism and } \kernel(v) \in \sf{Y}\}.
\end{align*}

Inspired by the work \cite{TRchen} of Chen and Lu,
for the $R$-bimodule $M$,
we consider the following Tor-vanishing condition:
\[\tag{$\sfT$} \Tor_{\geqslant 1}^R(M,M^{\otimes_Ri}\otimes_RP)=0
\text{ for each } P\in \Proj(R) \text{ and } i \geqslant 1,\]
which will play a key role in the proofs of the main results of the paper.
In this case, under some additional conditions on $M$,
we characterize Gorenstein projective (resp., injective, flat) modules over $T_R(M)$
in terms of the Gorenstein projective (resp., injective, flat) property of $R$-modules; see Corollaries \ref{G-proj} and \ref{G-inj} and Theorem \ref{G-flat}.
These results combine to the following main theorem of this paper.

\begin{intthm} \label{THM GP}
Let $M$ be an $N$-nilpotent $R$-bimodule.
\begin{rqm}
\item If $M$ satisfies the condition $(\sfT)$, $\pd_{R}M<\infty$ and $\fd_{\Rop}M<\infty$, then there is an equality $\GP(T_R(M))=\Phi(\GP(R))$.
\item If $M$ satisfies the condition $(\sfT)$, $\fd_{R}M<\infty$ and $\pd_{\Rop}M<\infty$, then there is an equality $\GI(T_R(M)^{\op})=\Psi(\GI(\Rop))$.
\item If $R$ is right coherent, $M$ is flat as an $R$-module and finitely presented as an $\Rop$-module, and $\pd_{\Rop}M<\infty$, then there is an equality $\GF(T_R(M))=\Phi(\GF(R))$.
\end{rqm}
\end{intthm}

Example \ref{example of qviver} shows that there exists a bimodule $M$ over some
algebra satisfying all the conditions in the above theorem.
We mention that the proof of Theorem \ref{THM GP}(3) relies on a characterization for the coherence of $T_R(M)$, which may be of independent interest to the reader; see Theorem \ref{coherence}.

\begin{intthm} \label{THE coherent}
Let $M$ be an $N$-nilpotent $R$-bimodule. Consider the next statements:
\begin{eqc}
\item $M$ is finitely generated as an $R$-module and $T_R(M)$ is left coherent.
\item $M$ is finitely presented as an $R$-module and $R$ is left coherent.
\end{eqc}
Then \eqclbl{i}$\!\implies\!$\eqclbl{ii} holds true. Furthermore, if $M$ is flat as an $\Rop$-module, then the two statements are equivalent.
\end{intthm}

We finish the paper with two applications of the above theorems.
We not only reobtain the earlier results in this direction,
but also get some new conclusions.
We first study Gorenstein projective (resp., injective, flat)
modules over trivial extension of rings, and get some new characterizations.
We then study when a module over a Morita context ring is Gorenstein projective (resp., injective, flat).
The corresponding results present some new characterizations of Gorenstein homological modules
over Morita context rings.

\section{Preliminaries}
\label{Preliminaries}
\noindent
In this section, we fix some notation, recall relevant notions
and collect some necessary facts.
Throughout the paper, we denote by $\Mod(R)$ the category of left $R$-modules,
and by $\Proj(R)$ (resp., $\Inj(R)$, $\Flat(R)$) the subcategory of $\Mod(R)$
consisting of all projective (resp., injective, flat) $R$-modules. For an $R$-module $X$, denote by $\pd_{R}X$ (resp., $\id_{R}X$, $\fd_{R}X$)
the projective (resp., injective, flat) dimension of $X$, and
by $X^+$ the character module $\Hom_\mathbb{Z}(X, \mathbb{Q/Z})$.

\begin{bfhpg}[\bf Tensor rings]\label{Tensor rings}
Let $R$ be a ring and $M$ an $R$-bimodule.
We write $M^{\otimes_R 0} = R$ and $M^{\otimes_R (i+1)} = M \otimes_R (M^{\otimes_R i})$
for $i\geqslant 0$.
For an integer $N \geqslant 0$, recall that $M$ is said to be $N$-\emph{nilpotent} if
$M^{\otimes_R(N+1)}=0$.
For an $N$-nilpotent $R$-bimodule $M$, we denote by $T_R(M)=\bigoplus_{i=0}^N M^{\otimes_Ri}$
the \emph{tensor ring} with respect to $M$. It is easy to check that $T_R(M)^{\op}\is T_{\Rop}(M)$.

Let $\Gamma$ be the category whose objects are the pairs $(X, u)$,
where $X \in \Mod(R)$ and $u: M\otimes_R X \to X$ is an $R$-homomorphism; the morphisms from $(X, u)$ to $(X', u')$ are
those $R$-homomorphisms $f \in \Hom_R(X, X')$ such that
$f \circ u = u' \circ (M\otimes f)$.
It follows from \cite{TRchen} that the category $\Mod(T_R(M))$ of $T_R(M)$-modules
is equivalent to the category $\Gamma$.

Dually, Let $\Omega$ be the category whose objects are the pairs $[Y, v]$,
where $Y \in \Mod(\Rop)$ and $v: Y \to \Hom_{\Rop}(M, Y)$ is an $\Rop$-homomorphism; the morphisms from $[Y, v]$ to $[Y', v']$ are
those $\Rop$-homomorphisms $g \in \Hom_{\Rop}(Y, Y')$
such that $\Hom_{\Rop}(M, g) \circ v = v' \circ g$. Then the category $\Mod(T_R(M)^{\op})$ of $T_R(M)^{\op}$-modules is equivalent to the category $\Omega$.

Unless otherwise specified, in the paper, we always view a  $T_R(M)$-module as a pair $(X,u)$ with $X \in \Mod(R)$ and $u \in \Hom_R(M \otimes_R X, X)$,
and view a $T_R(M)^{\op}$-module as a pair $[Y, v]$
with $Y \in \Mod(\Rop)$ and $v \in \Hom_{\Rop}(Y,\Hom_{\Rop}(M,Y))$. Note that a $T_R(M)^{\op}$-module can be equivalently viewed as a pair $(Y,\overline{v})$,
where $\overline{v} \in \Hom_{\Rop}(Y\otimes_RM, Y)$ is the adjoint morphism of $v$.

We mention that a sequence
$$(X, u) \overset{f} \longrightarrow (X', u') \overset{g} \longrightarrow (X'', u'')$$
in $\Mod(T_R(M))$  is exact if and only if
the underlying sequence
$X \overset{f} \longrightarrow X' \overset{g} \longrightarrow X''$
is exact in $\Mod(R)$.
\end{bfhpg}

\setup Throughout this paper, we always let $M$ denote an $N$-nilpotent $R$-bimodule.

\begin{bfhpg}[\bf The forgetful functor and its adjoint]
\label{The forgetful functor and its}
There exists a \emph{forgetful functor}
\begin{center}
$U: \Mod(T_R(M))\to \Mod(R)$,
\end{center}
which maps a  $T_R(M)$-module $(X,u)$ to the underlying $R$-module $X$.
Recall from \lemcite[2.1]{TRchen} that $U$ admits a left adjoint
\begin{center}
$\Ind: \Mod(R) \to \Mod(T_R(M))$,
\end{center}
defined as follows:
\begin{itemize}
\item For an $R$-module $X$, define
$\Ind(X)=(\bigoplus_{i=0}^N(M^{\otimes_Ri}\otimes_RX),c_X)$,
where $c_X$ is an inclusion from
$M\otimes_R(\bigoplus_{i=0}^N(M^{\otimes_Ri}\otimes_RX))\cong
\bigoplus_{i=1}^N(M^{\otimes_Ri}\otimes_RX)$
to
$\bigoplus_{i=0}^N(M^{\otimes_Ri}\otimes_RX)$.
More explicitly,
\begin{center}
$c_X=\begin{pmatrix}
0 & 0 & 0 & \cdots & 0 \\
1 & 0 & 0 & \cdots & 0 \\
0 & 1 & 0 & \cdots & 0 \\
\cdots & \cdots & \cdots & \cdots & \cdots \\
0 & 0 & 0 & \cdots & 1 \\
\end{pmatrix}_{(N+1)\times N}.$
\end{center}
\item For an $R$-homomorphism $f: X \to Y$,
the $T_R(M)$-homomorphism $\Ind(f): \Ind(X) \to \Ind(Y)$ can be viewed as
a formal diagonal matrix with diagonal elements
$M^{\otimes_Ri}\otimes f$ for $0 \leqslant i\leqslant N$.
\end{itemize}
\end{bfhpg}

\begin{bfhpg}[\bf The stalk functor and its adjoint]
\label{The stalk functor and its adjoints}
There exists a \emph{stalk functor}
\begin{center}
$S: \Mod(R)\to \Mod(T_R(M))$,
\end{center}
which maps an $R$-module $X$ to the  $T_R(M)$-module $(X, 0)$.
The functor $S$ admits a left adjoint
$C: \Mod(T_R(M)) \to \Mod(R)$ defined as follows:
\begin{itemize}
\item For a  $T_R(M)$-module $(X,u)$, define $C((X,u))=\cok(u)$.
\item For a morphism $f: (X,u) \to (Y,v)$ in $\Mod(T_R(M))$, the morphism
$C(f) : \cok(u) \to \cok(v)$ is induced by the universal property of cokernels.
\end{itemize}
\end{bfhpg}

\begin{rmk}\label{adjoint pairs}
By \ref{The forgetful functor and its} and \ref{The stalk functor and its adjoints},
there are adjoint pairs $(C, S)$ and $(\Ind, U)$ as follows:
  \begin{equation*}
  \xymatrix@C=3pc{
    \Mod(R)
    \ar[r]^-{S}
    &
    \Mod(T_R(M))
    \ar@/_1.8pc/[l]_-{C}
    \ar[r]^-{U}
    &
    \Mod(R)
    \ar@/_1.8pc/[l]_-{\Ind}.
  }
  \end{equation*}
It is easy to see that $C\circ\Ind=\Id_{\Mod(R)}$, and
it follows from the Eilenberg-Watts theorem that the functor $\Ind$ is isomorphic to the functor $T_R(M)\otimes_R-$, so one has $\Ind(X)\is T_R(M)\otimes_RX$ for each $R$-module $X$.

For a $T_R(M)$-module $(X,u)$,
there exists a short exact sequence
\[
\xymatrix{
  0 \to \Ind(M\otimes_RX)
   \ar[r]^-{\phi_{(X,u)}}
& \Ind(X)
   \ar[r]^-{\varepsilon_{(X,u)}}
& (X,u) \to 0}
\]
of $T_R(M)$-modules, where
\[
\phi_{(X,u)}=\begin{pmatrix}
-u & 0 & 0 & \cdots & 0\\
1 & -M\otimes u & 0 & \cdots & 0\\
0 & 1 & -M^{\otimes_R2}\otimes u & \cdots & 0\\
\cdots & \cdots & \cdots & \cdots & \cdots\\
0 & 0 & 0 & \cdots & 1\\
\end{pmatrix}_{(N+1)\times N}
\]
and
$$\varepsilon_{(X,u)}=\begin{pmatrix}
1 & u & u\circ (M\otimes u) & \cdots  & u\circ (M\otimes u)\circ \cdots \circ
(M^{\otimes_R(N-1)}\otimes u)
\end{pmatrix}$$
is the component of the counit $\varepsilon$ of the adjoin pair $(\Ind, U)$ at $(X,u)$; see \pgcite[191]{TRchen}.
\end{rmk}

\begin{ipg}\label{adjoint pairs-dual}
Similarly to \ref{The forgetful functor and its} and \ref{The stalk functor and its adjoints}, there are adjoint pairs $(S, K)$ and $(U, \Coind)$ as follows:
  \begin{equation*}
  \xymatrix@C=3pc{
    \Mod(\Rop)
    \ar[r]^-{S}
    &
    \Mod(T_R(M)^{\op})
    \ar@/^1.8pc/[l]_-{K}
    \ar[r]^-{U}
    &
    \Mod(\Rop).
    \ar@/^1.8pc/[l]_-{\Coind}
  }
  \end{equation*}
Here, $K([Y,v])=\kernel(v)$ for a $T_R(M)^{\op}$-module $[Y,v]$, and
$$\Coind(Y)=[\bigoplus_{i=0}^N\Hom_{\Rop}(M^{\otimes_Ri},Y),r_Y]$$
for an $\Rop$-module $Y$, where $r_Y$ is the morphism from
$\bigoplus_{i=0}^N\Hom_{\Rop}(M^{\otimes_Ri},Y)$ to
\[\Hom_{\Rop}(M,\bigoplus_{i=0}^N\Hom_{\Rop}(M^{\otimes_Ri},Y)) \cong
\bigoplus_{i=1}^N\Hom_{\Rop}(M^{\otimes_Ri},Y)\]
defined as
\[
r_Y=\begin{pmatrix}
0 & 1 & 0 & 0 & \cdots & 0 \\
0 & 0 & 1 & 0 & \cdots & 0 \\
0 & 0 & 0 & 1 & \cdots & 0 \\
\cdots & \cdots & \cdots & \cdots & \cdots & \cdots \\
0 & 0 & 0 & 0 & \cdots & 1 \\
\end{pmatrix}_{N\times(N+1)}.
\]
It is easy to see that $K\circ\Coind=\Id_{\Mod(\Rop)}$, and
it follows from the Eilenberg-Watts theorem that the functor $\Coind$ is isomorphic to the functor $\Hom_{\Rop}(T_R(M),-)$, so one has $\Coind(Y)\is \Hom_{\Rop}(T_R(M),Y)$ for each $\Rop$-module $Y$.

For a $T_R(M)^{\op}$-module $[Y, v]$, there exists a short exact sequence
\[
\xymatrix{
   0 \to [Y,v]
   \ar[r]^-{\eta_{[Y,v]}}
&  \Coind(Y)
   \ar[r]^-{\psi_{[Y,v]}}
&  \Coind(\Hom_{\Rop}(M,Y)) \to 0}
\]
of $T_R(M)^{\op}$-modules, where
\[
\psi_{[Y,v]}=\begin{pmatrix}
-v & 1 & 0 & \cdots  & 0\\
0 & -\Hom_{\Rop}(M,v) & 1 & \cdots  & 0\\
0 & 0 & -\Hom_{\Rop}(M^{\otimes_R2},v) & \cdots  & 0\\
\cdots & \cdots & \cdots & \cdots & \cdots\\
0 & 0 & 0 & \cdots & 1\\
\end{pmatrix}_{N\times(N+1)},\]
and
\[
\eta_{[Y,v]}=(
1 \,\, v \,\, \Hom_{\Rop}(M,v)\circ v \,\, \cdots \,\,
\Hom_{\Rop}(M^{\otimes_R(N-1)},v)\circ \cdots \circ \Hom_{\Rop}(M,v) \circ v)^\tr
\]
is the component of the unit $\eta$ of the adjoin pair $(U, \Coind)$ at $[Y,v]$.
\end{ipg}

We collect some useful facts in the following lemma,
which will be used frequently in the sequel.

\begin{lem} \label{projectives in F-Rep}
The following statements hold.
\begin{prt}
\item A $T_R(M)$-module $(X,u)$ is finitely generated if and only if
$\cok(u)$ is a finitely generated $R$-module.
\item A $T_R(M)$-module $(X,u)$ is free if and only if
$(X,u) \cong \Ind(F)$ for some free $R$-module $F$.
\item A $T_R(M)$-module $(X,u)$ is finitely generated and free
if and only if $(X,u) \cong \Ind(F)$
for some finitely generated and free $R$-module $F$.
\item A $T_R(M)^{\op}$-module $[Y,v]=0$ if and only if $\kernel(v)=0$.
\item For each $T_R(M)^{\op}$-module $[Y,v]$,
the morphism $S\circ K([Y,v]) \to [Y,v]$ is an essential monomorphism.
\end{prt}
\end{lem}
\begin{prf*}
The statement (a) is proved similarly as in \prpcite[3.2.1]{GFMODB}.
The statement (b) holds since the $T_R(M)$-module $(X,u)$ is free if and only if
$(X,u) \cong T_R(M)^{(I)}\cong \Ind(R)^{(I)} \cong \Ind(R^{(I)})$ for some index set $I$;
see Remark \ref{adjoint pairs}.
Then the statement (c) is an immediate consequence of the statements (a) and (b).

(d) If $[Y,v]=0$, then it is clear that $\kernel(v)=0$.
Conversely, if $\kernel(v)=0$, then $v$ is a monomorphism, and hence
$$\Hom_{\Rop}(M^{\otimes_Ri},v): \Hom_{\Rop}(M^{\otimes_Ri},Y) \to \Hom_{\Rop}(M^{\otimes_R(i+1)},Y)$$
is a monomorphism for each $1\leqslant i \leqslant N$.
However, $\Hom_{\Rop}(M^{\otimes_R(N+1)},Y)=0$ as $M$ is $N$-nilpotent.
This yields that $\Hom_{\Rop}(M^{\otimes_RN},Y)=0$. Iteratively, one gets that $Y=0$.

(e) Note that there exists a commutative diagram
\[\xymatrix{
  0 \ar[r]^{} & \kernel(v) \ar[d]_{0} \ar[r]^{f} & Y  \ar[d]_{v}  \\
  0 \ar[r]^{} & \Hom_{\Rop}(M,\kernel(v)) \ar[r]^{\quad (M,f) \,\,} & \Hom_{\Rop}(M,Y).} \]
This induces the monomorphism $$f: S\circ K([Y,v])=[\kernel(v),0] \to [Y,v].$$
It remains to show that the submodule $[\kernel(v),0]$ of $[Y,v]$ is essential.
Indeed, let $[X,u]$ be a nonzero submodule of $[Y,v]$ with the inclusion $g$.
Note that there exists the monomorphism $f_1: [\kernel(u),0] \to [X,u]$.
Since $v \circ g \circ f_1=\Hom_{\Rop}(M,g) \circ u \circ f_1=0$,
it follows from the universal property of kernels that
there exists a $\Rop$-homomorphism $g_1: \kernel(u) \to \kernel(v)$ such that $f \circ g_1=g \circ f_1$.
This implies that the diagram
\[\xymatrix{ [\kernel(u),0] \ar[d]_{g_1} \ar[r]^{\quad f_1} & [X,u]  \ar[d]_{g}  \\
  [\kernel(v),0] \ar[r]^{\quad f} & [Y,v]} \]
in $\Mod(T_R(M)^{\op})$ is commutative.
To complete the proof, it suffices to show that this diagram is a pullback.
If we are done, then $[\kernel(u),0] \cong [X,u] \cap [\kernel(v),0]$
as both $f$ and $g$ are monomorphisms.
Note that $[\kernel(u),0]\neq 0$ by (d).
This yields that $[X,u] \cap [\kernel(v),0] \neq 0$, as desired.
Indeed, suppose that
\[\xymatrix{ [Z,w] \ar[d]_{g_2} \ar[r]^{f_2} & [X,u]  \ar[d]_{g}  \\
  [\kernel(v),0] \ar[r]^{\quad f} & [Y,v]} \]
is a commutative diagram in $\Mod(T_R(M)^{\op})$.
Note that $\Hom_{\Rop}(M,g) \circ u \circ f_2=v \circ g \circ f_2=v \circ f \circ g_2=0$.
Since $\Hom_{\Rop}(M,g)$ is a monomorphism, one gets that $u \circ f_2=0$.
This implies that there exists a unique
$\Rop$-homomorphism $h: Z \to \kernel(u)$ such that $f_1 \circ h=f_2$.
By the fact that both $\Hom_{\Rop}(M,f_1)$ and $f$ are monomorphisms,
one can similarly verify that $\Hom_{\Rop}(M,h) \circ w=0$ and $g_1 \circ h=g_2$, respectively.
Thus, $h: [Z,w] \to [\kernel(u),0]$ is the desired $T_R(M)^{\op}$-homomorphism.
\end{prf*}

We mention that the following definition of the subcategory $\Phi(\sf{X})$ is essentially inspired by the classical works \cite{Bir35} by Birkhoff and \cite{AR91} by Auslander and Reiten.
Zhang \cite{2011ZHANG} refers to $\Phi({\sf{X}})$ as a monomorphism category.

\begin{ipg}
Let $\sf{X}$ be a subcategory of $R$-modules
and $\sf{Y}$ a subcategory of $R^{\op}$-modules.
The symbol $\Phi(\sf{X})$ (resp., $\Psi(\sf{Y})$) denotes the subcategory of $\Mod(T_R(M))$ (resp., $\Mod(T_R(M)^{\op})$) defined as follows:
\begin{align*}
&\Phi({\sf{X}}) = \{(X,u) \in \Mod(T_R(M)) \,|\,
u \textrm{ is a monomorphism and } \cok(u) \in \sf{X}\}; \\
&\Psi({\sf{Y}}) = \{[Y,v] \in \Mod(T_R(M)^{\op}) \,|\,
v \textrm{ is an epimorphism and } \kernel(v) \in \sf{Y}\}.
\end{align*}
\end{ipg}

\begin{lem} \label{homological modules}
The following equalities hold.
\begin{prt}
\item $\Proj(T_R(M))=\Ind(\Proj(R))=\Phi(\Proj(R))$.
\item $\Inj(T_R(M)^{\op})=\Coind(\Inj(R^{\op}))=\Psi(\Inj(R^{\op}))$.
\item $\Flat(T_R(M))=\Phi(\Flat(R))$.
\end{prt}
\end{lem}
\begin{prf*}
(a) The equality $\Ind(\Proj(R))=\Phi(\Proj(R))$ holds by \lemcite[2.3]{TRchen},
and it is clear that $\Ind(\Proj(R))\subseteq\Proj(T_R(M))$ as $(\Ind,U)$ is an adjoint pair with $U$ exact; see \ref{The forgetful functor and its}.
Thus, we have to prove the inclusion $\Proj(T_R(M)) \subseteq \Phi(\Proj(R))$.

Let $(X,u)$ be a projective $T_R(M)$-module. Then one gets that $(X,u)\oplus(Y,v)\cong \Ind(R^{(I)})$ for some $T_R(M)$-module $(Y,v)$ and index set $I$
by Lemma \ref{projectives in F-Rep}(b),
and so there exists a commutative diagram with split exact rows:
\[
\xymatrix{
  0 \ar[r]^{}
  & M\otimes_RX  \ar[d]_{u}\ar[r]^{}
  & \bigoplus_{i=1}^N(M^{\otimes_Ri}\otimes_RR^{(I)})
    \ar[d]_{c_{R^{(I)}}}  \ar[r]^{}
  & M\otimes_RY \ar[d]_{v}
    \ar[r]^{}  & 0  \\
  0 \ar[r]^{}
  & X   \ar[r]^{}
  & \bigoplus_{i=0}^N(M^{\otimes_Ri}\otimes_RR^{(I)})
    \ar[r]^{} & Y  \ar[r]^{}  & 0. }
\]
This implies that $u$ is a monomorphism as so is $c_{R^{(I)}}$. We mention that $\cok(u)=C(X,u)$ is a projective $R$-module as $(C,S)$ is an adjoint pair with $S$ exact; see \ref{The stalk functor and its adjoints}.
Thus, one has $(X,u) \in \Phi(\Proj(R))$.

(b) By a dual argument of the proof of \lemcite[2.3]{TRchen},
one sees that the equality $\Coind(\Inj(R^{\op}))=\Psi(\Inj(R^{\op}))$ holds.
Note that the functor $\Coind$ preserves injective modules; see \ref{adjoint pairs-dual}.
This implies that $\Coind(\Inj(R^{\op})) \subseteq \Inj(T_R(M)^{\op})$.
Thus, we only need to prove the inclusion
$\Inj(T_R(M)^{\op}) \subseteq \Coind(\Inj(R^{\op}))$.
Indeed, let $[Y,v]$ be an injective $T_R(M)^{\op}$-module.
Since the functor $K$ also preserves injective modules, the $R^{\op}$-module $\kernel(v)$ is injective.
By Lemma \ref{projectives in F-Rep}(e), we see that
both $[Y,v]$ and $\Coind(\kernel(v))$ are injective envelopes of $[\kernel(v),0]$,
and hence, $[Y,v] \cong \Coind(\kernel(v))\in \Coind(\Inj(R^{\op}))$, as desired.

(c) It is easy to check that a $T_R(M)$-module $(X,u)$ is in $\Phi(\Flat(R))$ if and only if the character module $(X,u)^+=[X^+,u^+]$ is in $\Psi(\Inj(\Rop))$.
Thus, the equality holds by the equalities in (b).
\end{prf*}

The following results are immediate consequences of Lemma \ref{homological modules}.

\begin{cor} \label{Proj of Ind}
Let $X$ be an $R$-module. Then the following statements hold.
\begin{prt}
\item $X$ is projective if and only if $\Ind(X)$ is a projective $T_R(M)$-module.
\item If $\Tor_{\geqslant 1}^R(T_R(M),X)=0$, then $\pd_RX=\pd_{T_R(M)}\Ind(X)$.
\end{prt}
\end{cor}

\begin{cor} \label{Inj of Ind-}
Let $Y$ be an $\Rop$-module. Then the following statements hold.
\begin{prt}
\item $Y$ is injective if and only if $\Coind(Y)$ is an injective $T_R(M)^{\op}$-module.
\item If $\Ext^{\geqslant 1}_{\Rop}(T_R(M),Y)=0$, then $\id_{\Rop}Y=\id_{T_R(M)^{\op}}\Coind(Y)$.
\end{prt}
\end{cor}

\section{Gorenstein projective and injective $T_R(M)$-modules}
\label{Gorenstein projective  T_R(M)-modules}
\noindent
In this section, we give the proofs of (1) and (2) in Theorem \ref{THM GP} in the introduction. We begin with the following definitions.

\begin{bfhpg}[\bf Gorenstein projective and injective modules]
\label{Gorenstein homological modules}
Recall from \dfncite[2.1]{GHD} that an exact sequence
$$\xymatrix@C=0.5cm{
P^\bullet: \cdots \to P^{-1} \ar[r]^{\qquad \quad d^{-1} \,} & P^0
\ar[r]^{\,\, d^0 \qquad } & P^1 \to \cdots,}$$
of projective  $R$-modules is called \emph{complete projective resolution} if it remains exact after applying the functor $\Hom_R(-,P)$ for every projective $R$-module $P$.
An $R$-module $X$ is called \emph{Gorenstein projective}
provided that there exists a complete projective resolution $P^\bullet$ such that
$X \cong \kernel(P^0 \to P^1)$.
Dually, an exact sequence
$$\xymatrix@C=0.5cm{
I^\bullet: \cdots \to I^{-1} \ar[r]^{\qquad \quad d^{-1} \,\,} & I^0
\ar[r]^{\,\, d^0 \qquad } & I^1 \to \cdots,}$$
of injective $R$-modules is called
\emph{complete injective resolution} if it remains exact after applying
the functor $\Hom_{R}(E,-)$ for every injective $R$-module $E$.
An $R$-module $X$ is called \emph{Gorenstein injective}
provided that there exists a complete injective resolution $I^\bullet$
such that $X \cong \kernel(I^0 \to I^1)$. We denote by $\GP(R)$ and $\GI(R)$ the subcategory of Gorenstein projective and Gorenstein injective $R$-modules, respectively.
\end{bfhpg}

\begin{dfn} \label{dfn of compatible}
An $R$-bimodule $X$ is said to be \emph{compatible}
if the following two conditions hold
for a complete projective resolution $Q^\bullet$ of $T_R(M)$-modules:
\begin{prt}
\item[(C1)]
The complex $X\otimes_RQ^\bullet$ is exact.
\item[(C2)]
The complex $\Hom_{T_R(M)}(Q^\bullet,\Ind(X\otimes_RP))$ is exact
for each $P \in \Proj(R)$.
\end{prt}
\end{dfn}

\begin{rmk}\label{rmk of compatible}
Compatible bimodules were first defined in \dfncite[1.1]{2013ZHANG} to describe
Gorenstein projective modules over triangular matrix Artin algebras.
This notion were further extened in \cite{GFGN} for
Morita context Artin algebras with zero bimodule homomorphisms.
Recently, Guo and Xi \cite{2022XI} weakened the compatibility condition
to construct Gorenstein projective modules over
Morita context rings with one bimodule homomorphism zero.
Moreover, the notion of generalized compatible bimodules was introduced in \dfncite[3.1]{TriExtGPMao} to describe Gorenstein projective modules over trivial ring extensions. Our Definition \ref{dfn of compatible} is inspired by \prpcite[3.8]{TRchen}.
\end{rmk}

The following definition essentially come from \dfncite[3.5]{TRchen}.

\begin{dfn}
An $R$-bimodule $X$ is said to be \emph{admissible} if
$$\Ext_R^1(G,X^{\otimes_Ri}\otimes_RP)=0=\Tor^R_1(X, X^{\otimes_Ri}\otimes_RG)$$
for each $G \in \GP(R)$, $P \in \Proj(R)$ and $i\geqslant 0$.
\end{dfn}

Recall that the condition $(\sfT)$ for an $R$-bimodule $X$ is defined as follows:
\[\tag{$\sfT$} \Tor_{\geqslant 1}^R(X,X^{\otimes_Ri}\otimes_RP)=0
\text{ for each } P\in \Proj(R) \text{ and } i \geqslant 1.\]
In the following, we prove that if the $R$-bimodule $M$ with $\pd_R M<\infty$ and $\fd_{\Rop}M<\infty$ satisfies the condition $(\sfT)$,
then it is compatible and admissible.
We need the following two auxiliary results.

The next result can be proved similarly as in \lemcite[4.2]{TRchen} (see \lemcite[6.5]{KP24}).

\begin{lem} \label{lemma4.2 for (T)}
Let $X$ be an $R$-bimodule satisfying the condition $(\sfT)$.
Then the following statements are equivalent for each $R$-module $Y$.
\begin{eqc}
\item
$\Tor_{\geqslant 1}^R(X,X^{\otimes_Ri}\otimes_RY)=0$ for any $i \geqslant 0$.
\item
$\Tor_{\geqslant 1}^R(X^{\otimes_Rs},X^{\otimes_Ri}\otimes_RY)=0$
for any $s \geqslant 1$ and $i \geqslant 0$.
\item
$\Tor_{\geqslant 1}^R(X^{\otimes_Rs},Y)=0$ for any $s \geqslant 1$.
\end{eqc}
\end{lem}

One can prove the following result by using a similar method
as in the proof of \lemcite[4.5]{TRchen}.

\begin{lem} \label{lemma4.5 for (T)}
Let $X$ be an $R$-bimodule satisfying the condition $(\sfT)$.
If $\pd_RX$ $($resp., $\pd_{\Rop}X$, $\fd_RX$, $\fd_{\Rop}X$$)$ is finite,
then so is $\pd_RX^{\otimes_Ri}$ $($resp., $\pd_{\Rop}X^{\otimes_Ri}$,
$\fd_RX^{\otimes_Ri}$, $\fd_{\Rop}X^{\otimes_Ri}$$)$
for each $i \geqslant 1$.
\end{lem}

\begin{prp} \label{example of admissible and compatible}
Suppose that the $R$-bimodule $M$ satisfies the condition $(\sfT)$.
If $\pd_RM<\infty$ and $\fd_{\Rop}M<\infty$,
then $M$ is compatible and admissible.
\end{prp}
\begin{prf*}
Let $Q^\bullet$ be a complete projective resolution of projective $T_R(M)$-modules.
By Lemma \ref{homological modules}(a), $Q^\bullet$ can be written as:
$$\xymatrix@C=0.5cm{
\cdots \to \Ind(P^{-1}) \ar[r]^{\qquad d^{-1}\,\,} & \Ind(P^0)
\ar[r]^{d^0 \quad \,\,} & \Ind(P^1) \to \cdots,}$$
where $P^j \in \Proj(R)$ for each $j \in \mathbb{Z}$. Hence, the complex
$$P^\bullet:\ \xymatrix@C=0.5cm{\cdots \to \bigoplus_{i=0}^N(M^{\otimes_Ri}\otimes_RP^{-1})
\ar[r]^{\,\, d^{-1}} &
\bigoplus_{i=0}^N(M^{\otimes_Ri}\otimes_RP^0) \to \cdots}$$ of $R$-modules is exact.
To prove the condition (C1), we have to prove that $M\otimes_RP^\bullet$ is exact.
Fix a flat resolution
$0 \to F^n \to \cdots \to  F^1 \to  F^0 \to M \to 0$ of $\Rop$-module $M$.
We mention that $\Tor_{\geqslant 1}^R(M,\bigoplus_{i=0}^N(M^{\otimes_Ri}\otimes_RP^j))=0$
as $M$ satisfies the condition $(\sfT)$.
Then there exists an exact sequence of complexes
$$0 \to F^n\otimes_RP^\bullet \to \cdots \to  F^1\otimes_RP^\bullet \to
F^0\otimes_RP^\bullet \to M\otimes_RP^\bullet \to 0.$$
Since each $F^k\otimes_RP^\bullet$ is exact for $0\leqslant k\leqslant n$,
one gets that $M\otimes_RP^\bullet$ is exact, as desired.

We then prove the condition (C2). It follows from Lemma \ref{lemma4.2 for (T)} that
$$\Tor_{\geqslant 1}^R(M^{\otimes_Ri},M\otimes_RP)=0$$
for each $P \in \Proj(R)$ and $1\leqslant i\leqslant N$,
and hence, $\Tor_{\geqslant 1}^R(T_R(M),M\otimes_RP)=0$.
Therefore, by Corollary \ref{Proj of Ind}(b), one has $\pd_{T_R(M)}\Ind(M\otimes_RP)=\pd_R(M\otimes_RP)<\infty$ as $\pd_RM<\infty$. This yields that $\Hom_{T_R(M)}(Q^\bullet,\Ind(M\otimes_RP))$ is exact.
Thus, $M$ is compatible.

Next, we prove that $M$ is admissible.
Indeed, by Lemma \ref{lemma4.5 for (T)}, one gets that $\pd_R(M^{\otimes_Rs}\otimes_RP)<\infty$ for each $0\leqslant s \leqslant N$ as $\pd_RM<\infty$. This yields that $\Ext_R^{\geqslant 1}(G,M^{\otimes_Rs}\otimes_RP)=0$
for each $G \in \GP(R)$.
On the other hand, one has $\fd_{\Rop}(M^{\otimes_Rs})<\infty$ by Lemma \ref{lemma4.5 for (T)},
which ensures that $\Tor_{\geqslant 1}^R(M^{\otimes_Rs}, G)=0$.
Thus, by Lemma \ref{lemma4.2 for (T)},
one gets that $\Tor_{\geqslant 1}^R(M,M^{\otimes_Rs}\otimes_RG)=0$. This yields that $M$ is admissible.
\end{prf*}

In the following, we prove that if $M$ is compatible and admissible,
then one has $\GP(T_R(M))=\Phi(\GP(R))$.

\begin{lem} \label{GPTRM in Phi}
Suppose that the $R$-bimodule $M$ is compatible. Then there is an inclusion
$\GP(T_R(M)) \subseteq \Phi(\GP(R))$.
\end{lem}

\begin{prf*}
Let $(X,u)$ be in $\GP(T_R(M))$ with
$$Q^\bullet:\ \xymatrix@C=0.5cm{
\cdots \to \Ind(P^{-1}) \ar[r]^-{d^{-1}} & \Ind(P^{0})
\ar[r]^-{d^0} & \Ind(P^{1}) \to \cdots}$$
a complete projective resolution of projective $T_R(M)$-modules such that
$(X,u)\is\ker(d^0)$, where each $P^i$ is in $\Proj(R)$.
Then there is an exact complex
$$\xymatrix@C=0.5cm{
\cdots \to
\bigoplus_{i=0}^N(M^{\otimes_Ri}\otimes_RP^{-1})
\ar[r]^{\,\,d^{-1}} & \bigoplus_{i=0}^N(M^{\otimes_Ri}\otimes_RP^0)
\to \cdots}$$
in $\Mod(R)$.
Since $M\otimes_RQ^\bullet$ is exact by assumption, the complex
$$\xymatrix@C=1cm{\cdots \to \bigoplus_{i=1}^N(M^{\otimes_Ri}\otimes_RP^{-1})
\ar[r]^{\,\, M\otimes d^{-1}} &
\bigoplus_{i=1}^N(M^{\otimes_Ri}\otimes_RP^0) \to \cdots}$$
is exact. Thus, one gets a commutative diagram with exact rows and columns:
$$
\xymatrix@C=0.7cm@R=0.7cm{
  0 \ar[r]^{}
  & \bigoplus_{i=1}^N(M^{\otimes_Ri}\otimes_RP^{-1}) \ar[d]_{M\otimes d^{-1}}
  \ar[r]^{c_{P^{-1}}}
  & \bigoplus_{i=0}^N(M^{\otimes_Ri}\otimes_RP^{-1})
    \ar[d]_{d^{-1}} \ar[r]^{\qquad \qquad \pi_{P^{-1}}}
  & P^{-1} \ar[d]_{} \ar[r]^{} & 0  \\
  0 \ar[r]^{}
  & \bigoplus_{i=1}^N(M^{\otimes_Ri}\otimes_RP^0) \ar[d]_{M\otimes d^0}
    \ar[r]^{c_{P^0}}
  & \bigoplus_{i=0}^N(M^{\otimes_Ri}\otimes_RP^0)
    \ar[d]_{d^0} \ar[r]^{\qquad \qquad \pi_{P^0}}
  & P^0 \ar[d]_{} \ar[r]^{} & 0  \\
  0 \ar[r]^{}
  & \bigoplus_{i=1}^N(M^{\otimes_Ri}\otimes_RP^1)
  \ar[r]^{c_{P^1}}
  & \bigoplus_{i=0}^N(M^{\otimes_Ri}\otimes_RP^1)
  \ar[r]^{\qquad \qquad \pi_{P^1}}
  & P^1 \ar[r]^{} & 0. }$$
This induces a short exact sequence $0 \to M\otimes_RX \xra{u} X \to \cok(u) \to 0$.
Now we only need to prove that $\cok(u)\in\GP(R)$.
It is sufficient to show that $\Hom_R(P^\bullet,P)$ is exact for each $P \in \Proj(R)$,
where $P^\bullet$ is the third non-zero column in the above commutative diagram.
We mention that there exists a short exact sequence
$$\xymatrix@C=1cm{0 \to
\Ind(M\otimes_RP) \ar[r]^-{\phi_{(P,0)}} &
\Ind(P) \ar[r]^-{\eta_{(P,0)}} & (P,0) \to 0}$$
in $\Mod(T_R(M))$; see Remark \ref{adjoint pairs}.
Applying the functor $\Hom_{T_R(M)}(Q^\bullet, -)$ to this sequence, one gets a short exact sequence
\[\xymatrix@C=1.5cm{ 0 \longrightarrow  (Q^\bullet,\Ind(M\otimes_RP))\ar[r]^-{(Q^\bullet,\phi_{(P,0)})} &
(Q^\bullet,\Ind(P))
 \ar[r]^-{(Q^\bullet,\eta_{(P,0)})} &(Q^\bullet,(P,0)) \longrightarrow 0}
\]
of complexes as each term of $Q^\bullet$ is projective.
Note that the complex $Q^\bullet$ is a complete projective resolution,
so $\Hom_{T_R(M)}(Q^\bullet,\Ind(P))$ is exact.
By assumption, $\Hom_{T_R(M)}(Q^\bullet,\Ind(M\otimes_RP))$ is also exact.
Thus, $\Hom_{T_R(M)}(Q^\bullet, S(P))=\Hom_{T_R(M)}(Q^\bullet,(P,0))$ is exact,
which yields that $\Hom_R(C(Q^\bullet),P)$ is exact.
However, the complex $C(Q^\bullet)$ is nothing but $P^\bullet$
as $C\circ\Ind=\Id_{\Mod(R)}$; see Remark \ref{adjoint pairs}.
Cosequently, $\Hom_R(P^\bullet,P)$ is exact, as desired.
\end{prf*}

Before proving the other inclusion
$\Phi(\GP(R)) \subseteq \GP(T_R(M))$,
we finish a few preparatory works.

\begin{lem} \label{lem of Tor vanish}
Let $(X,u)$ be a $T_R(M)$-module such that $u$ is a monomorphism.
If $\Tor_1^R(M,M^{\otimes_Ri}\otimes_R\cok(u))=0$ for each $0\leqslant i \leqslant N$,
then $\Tor_1^R(M,X)=0$.
\end{lem}

\begin{prf*}
Applying the functor $M\otimes_R-$ to the exact sequence
$$\delta:\ 0 \to M\otimes_RX \overset{u\,\,}\longrightarrow X
\overset{\pi \,}\longrightarrow \cok(u) \to 0,$$
one gets exact sequences
$$M\otimes \delta:\ 0 \to M\otimes_R(M\otimes_RX) \to M\otimes_RX \to M\otimes_R\cok(u) \to 0$$
and
$$\lambda_1:\ \cdots \to \Tor^R_1(M, M\otimes_RX) \to \Tor^R_1(M,X) \to 0$$
as $\Tor^R_1(M,\cok(u))=0$ by assumption. Similarly, applying the functor $M\otimes_R-$ to the sequence $M\otimes \delta$, one gets exact sequences
$$M^{\otimes_R2}\otimes \delta:\ 0 \to M^{\otimes_R3}\otimes_RX \to
M^{\otimes_R2}\otimes_RX \to M^{\otimes_R2}\otimes_R\cok(u) \to 0$$
and
$$\lambda_2:\ \cdots \to \Tor^R_1(M, M^{\otimes_R2}\otimes_RX)
\to \Tor^R_1(M,M\otimes_RX) \to 0$$
as $\Tor^R_1(M, M\otimes_R\cok(u))=0$. In the same way, for each $3 \leqslant j \leqslant N$, one gets exact sequences
$$M^{\otimes_Rj}\otimes \delta:\ 0  \to M^{\otimes_R(j+1)}\otimes_RX \to
M^{\otimes_Rj}\otimes_RX \to M^{\otimes_Rj}\otimes_R\cok(u) \to 0$$
and
$$\lambda_j:\ \cdots \to \Tor^R_1(M, M^{\otimes_Rj}\otimes_RX)
\to \Tor^R_1(M,M^{\otimes_R(j-1)}\otimes_RX) \to 0.$$
Since $M$ is $N$-nilpotent, it follows that
$M^{\otimes_RN}\otimes_RX \cong M^{\otimes_RN}\otimes_R\cok(u)$.
Hence, one has
$\Tor^R_1(M,M^{\otimes_RN}\otimes_RX) \cong
\Tor^R_1(M,M^{\otimes_RN}\otimes_R\cok(u))=0$, and so by the exact sequence $\lambda_N$,
we conclude that $\Tor^R_1(M,M^{\otimes_R(N-1)}\otimes_RX)=0$.
Iteratively, by the exact sequences $\lambda_{N-1}, \cdots, \lambda_1$,
one gets that $\Tor^R_1(M,X)=0$.
\end{prf*}

Recall that for a morphism $f: (X,u) \to (Y,v)$ in $\Mod(T_R(M))$, the morphism
$C(f): \cok(u) \to \cok(v)$ is induced by the universal property; see \ref{The stalk functor and its adjoints}.

\begin{lem} \label{f is mono}
Let $(X,u)$ and $(Y,v)$ be $T_R(M)$-modules
with $M^{\otimes_R j}\otimes v$ a monomorphism for each $0\leqslant j \leqslant N-1$.
If $f: (X,u) \to (Y,v)$ is a $T_R(M)$-homomorphism such that
$M^{\otimes_R i}\otimes C(f)$ is a monomorphism for every $0\leqslant i \leqslant N$,
then $f$ is a monomorphism.
\end{lem}
\begin{prf*}
The $T_R(M)$-homomorphism $f: (X,u) \to (Y,v)$ induces the commutative diagram
\[
\xymatrix{
  M\otimes_RX  \ar[d]_{M\otimes f} \ar[r]^{\quad \,\, u}
  & X \ar[d]_{f}  \ar[r]^{}
  & \cok(u) \ar@{>->}[d]_{C(f)} \ar[r]^{}  & 0  \\
  M\otimes_RY   \ar@{>->}[r]^{\quad \,\, v}  & Y
    \ar[r]^{} & \cok(v)  \ar[r]^{}  & 0. }
\]
Applying the functor $M^{\otimes_Rj}\otimes_R-$ to the above diagram
for each $0\leqslant j \leqslant N-1$,
one gets the commutative diagram
\[
\xymatrix@C=1.2cm{
   M^{\otimes_R(j+1)}\otimes_RX  \ar[d]_{M^{\otimes_R(j+1)}\otimes f}
      \ar[r]^{\,\,\,\,\, M^{\otimes_Rj}\otimes u}
  & M^{\otimes_Rj}\otimes_RX \ar[d]_{M^{\otimes_Rj}\otimes f}  \ar[r]^{}
  & M^{\otimes_Rj}\otimes_R\cok(u) \ar@{>->}[d]_{M^{\otimes_Rj}\otimes C(f)} \ar[r]^{}  & 0  \\
   M^{\otimes_R(j+1)}\otimes_RY   \ar@{>->}[r]^{\,\,\,\,\, M^{\otimes_Rj}\otimes v}
   & M^{\otimes_Rj}\otimes_RY
    \ar[r]^{} & M^{\otimes_Rj}\otimes_R\cok(v)  \ar[r]^{}  & 0.}
\]
It suffices to show that $M^{\otimes_RN}\otimes f$ is a monomorphism.
If we have done, using the Snake Lemma, one can iteratively prove that $f$ is a monomorphism.
Indeed, apply the functor $M^{\otimes_RN}\otimes_R-$ to the first diagram.
Note that $M^{\otimes_R(N+1)}\otimes_RY=0=M^{\otimes_R(N+1)}\otimes_RX$.
This implies that $M^{\otimes_RN}\otimes f$ is a monomorphism as
so is $M^{\otimes_RN}\otimes C(f)$ by assumption.
\end{prf*}

We can now prove the following result.

\begin{prp} \label{Phi in GPTRM}
Suppose that the $R$-bimodule $M$ is admissible.
Then there is an inclusion $\Phi(\GP(R)) \subseteq \GP(T_R(M))$.
\end{prp}
\begin{prf*}
Let $(X,u)$ be in $\Phi(\GP(R))$.
We prove that $(X,u)$ is Gorenstein projective.
To this end, we prove the following statements:
\begin{rqm}
\item There is an exact sequence
      $\eta:\ 0 \to (X,u) \to \Ind(P) \to (Y,v) \to 0$ in $\Mod(T_R(M))$
      with $P \in \Proj(R)$ and $(Y,v) \in \Phi(\GP(R))$ such that the sequence $\Hom_{T_R(M)}(\eta,\Ind(Q))$ is exact
      for each $Q \in \Proj(R)$.
\item There is an exact sequence
      $\xi:\ 0 \to (Z,w) \to \Ind(P') \to (X,u) \to 0$ in $\Mod(T_R(M))$
      with $P' \in \Proj(R)$ and $(Z,w) \in \Phi(\GP(R))$ such that the sequence $\Hom_{T_R(M)}(\xi,\Ind(Q))$ is exact for each $Q \in \Proj(R)$.
\end{rqm}

\textbf{For (1)}:
Consider the short exact sequence
\[
\delta:\ 0 \to M\otimes_RX \overset{u \,}\longrightarrow X
\overset{\pi \,}\longrightarrow \cok(u) \to 0
\]
in $\Mod(R)$.
Since $\cok(u) \in \GP(R)$,
there exists a short exact sequence
\[
\epsilon:\ 0 \to \cok(u) \overset{f}\longrightarrow P \overset{g}\longrightarrow H \to 0
\]
in $\Mod(R)$ with $P \in \Proj(R)$ and $H \in \GP(R)$.
Set $a_0= f \circ \pi : X \to P$, and hence,
we obtain an $R$-homomorphism
$M\otimes a_0= M\otimes_RX \to M\otimes_RP$.
Applying the functor $\Hom_R(-,M\otimes_RP)$ to $\delta$.
Since $\Ext_R^1(\cok(u),M\otimes_RP)=0$ by assumption,
there exists an $R$-homomorphism $a_1: X \to M\otimes_RP$ such that
$a_1 \circ u= M\otimes a_0$.
We iterate the above argument to obtain
\[
a_j:X \to M^{\otimes_Rj}\otimes_RP\]
such that $a_j \circ u= M\otimes a_{j-1}$ holds for $2 \leqslant j \leqslant N$.
Set $\widetilde{a}=(a_0, a_1, \cdots , a_N)^\tr$,
which is an $R$-homomorphism from
$X$ to $\bigoplus_{i=0}^N(M^{\otimes_Ri}\otimes_RP)$.
It is routine to check that $\widetilde{a} \circ u = c_P \circ (M\otimes \widetilde{a})$.
This implies that $\widetilde{a}=(a_0, a_1, \cdots , a_N)^\tr$
forms a $T_R(M)$-homomorphism from $(X,u)$ to $\Ind(P)$.
We mention that $\pi' \circ \widetilde{a}=a_0=f \circ \pi$. Then the diagram
\begin{equation*}
\xymatrix@C=0.5cm@R=1cm{
  0 \ar[r] & M\otimes_RX \ar[d]_{M\otimes \widetilde{a}}\ar[r]^{u} & X \ar[d]_{\widetilde{a}}\ar[r]^{\pi} & \cok(u) \ar[d]_{f}\ar[r] & 0  \\
  0 \ar[r] & \bigoplus_{i=1}^N(M^{\otimes_Ri}\otimes_RP) \ar[r]^{c_P} & \bigoplus_{i=0}^N(M^{\otimes_Ri}\otimes_RP) \ar[r]^{\qquad \qquad \pi'} & P \ar[r] & 0}
\end{equation*}
is commutative, which yields that $C(\widetilde{a})=f$. It is easy to see that $M^{\otimes_Rk}\otimes f$ is a monomorphisms for each $0 \leqslant k \leqslant N$ as $\Tor^R_1(M, M^{\otimes_Ri}\otimes_RH)=0$ for all $i\geqslant 0$ by assumption. On the other hand, since $c_P$ is a split monomorphism,
one gets that $M^{\otimes_Rj}\otimes c_P$ is a monomorphism
for each $0 \leqslant j \leqslant N-1$.
Thus, it follows from Lemma \ref{f is mono} that $\widetilde{a}$ is a monomorphism, and so there is a short exact sequence
\[
\eta:\ 0 \to (X,u) \overset{\widetilde{a}\,}\rightarrow \Ind(P)\to (Y,v) \to 0
\]
in $\Mod(T_R(M))$.
Next, we prove that $(Y,v) \in \Phi(\GP(R))$.
Note that the sequence $\eta$ can be written as
the following commutative diagram with exact rows and columns:
\begin{equation*} \label{factorization 1}
\tag{\ref{Phi in GPTRM}.1}
\xymatrix@C=1cm@R=0.5cm{
  & 0 \ar[d]^{} & 0 \ar[d]^{}  \\
  & M\otimes_RX \ar[d]_{u} \ar[r]^{M\otimes \widetilde{a} \,\,\, \qquad}
  & \bigoplus_{i=1}^N(M^{\otimes_Ri}\otimes_RP) \ar[d]_{c_P} \ar[r]^{}
  & M\otimes_RY \ar[d]_{v} \ar[r]^{} & 0  \\
  0 \ar[r]^{}
  & X \ar[d]_{\pi} \ar[r]^-{\widetilde{a}}
  & \bigoplus_{i=0}^N(M^{\otimes_Ri}\otimes_RP)
    \ar[d]_{\pi'} \ar[r]^{}
  & Y \ar[d]_{} \ar[r]^{} & 0  \\
  0 \ar[r]^{}
  & \cok(u) \ar[r]^{\quad f} \ar[d]^{}
  & P \ar[r]^{g} \ar[d]^{}
  & H \ar[r]^{}  & 0 .  \\
  & 0  & 0   }
\end{equation*}
It then follows from the Snake Lemma that the sequence
$0 \to M\otimes_RY \overset{v \,}\rightarrow Y \to H \to 0$
is exact, and so $(Y,v) \in \Phi(\GP(R))$, as desired.

Next, we prove that $\Hom_{T_R(M)}(\eta,\Ind(Q))$ is a short exact sequence for each $Q \in \Proj(R)$. It suffices to prove that $\Hom_{T_R(M)}(\widetilde{a},\Ind(Q))$ is an epimorphism.
To this end, fix a $T_R(M)$-homomorphism
$\widetilde{l}=(l_0,l_1, \cdots, l_N)^\tr: (X,u) \to \Ind(Q)$
with $l_i: X \to M^{\otimes_Ri}\otimes_RQ$ for $0 \leqslant i \leqslant N$. Then one has $l_0\circ u=0$ and $l_i\circ u=M\otimes l_{i-1}$ for $1 \leqslant i \leqslant N$.
Set $d_0=C(\widetilde{l}): \cok(u) \to Q$. Then one has $d_0 \circ \pi=l_0$.
Apply $\Hom_R(-,Q)$ to the sequence $\epsilon$.
Since $\Ext_R^1(H, Q)=0$,
there exists $b_0': P \to Q$ such that $b_0' \circ f= d_0$.
Set $b_0=b_0' \circ \pi':\bigoplus_{i=0}^N(M^{\otimes_Ri}\otimes_RP) \to Q$.
It is easy to check that
\[
b_0 \circ \widetilde{a} = b_0' \circ \pi' \circ \widetilde{a}
 = b_0' \circ f \circ \pi = d_0 \circ \pi = l_0 \text{ \quad and \quad }
b_0 \circ c_P = b_0' \circ \pi' \circ c_P = 0.
\]
Applying $\Hom_R(-,M\otimes_RQ)$ to the second non-zero column in (\ref{factorization 1}), one gets
\begin{center}
$b_1':\bigoplus_{i=0}^N(M^{\otimes_Ri}\otimes_RP) \to M\otimes_RQ$
\end{center} such that
$b_1' \circ c_P=M\otimes b_0$. Consider the $R$-homomorphism $l_1-b_1'\circ \widetilde{a}: X\to M\otimes_RQ$.
Then one gets that
\begin{align*}
(l_1-b_1'\circ \widetilde{a}) \circ u&=l_1\circ u-b_1'\circ \widetilde{a}\circ u \\
&=M\otimes l_0-b_1'\circ (c_P\circ (M\otimes \widetilde{a})) \\
&=M\otimes l_0-(M\otimes b_0)\circ (M\otimes \widetilde{a}) \\
&=M\otimes l_0-M\otimes(b_0 \circ \widetilde{a}) \\
&=M\otimes l_0-M\otimes l_0 \\
&=0.
\end{align*}
By the universal property of cokernels,
there exists
$x_1: \cok(u) \to M\otimes_RQ$ such that $x_1\circ \pi=l_1-b_1'\circ \widetilde{a}$.

Apply $\Hom_R(-,M\otimes_RQ)$ to the sequence $\epsilon$.
Since $\Ext_R^1(H, M\otimes_RQ)=0$ by assumption,
there exists $y_1:P \to M\otimes_RQ$ such that $y_1 \circ f= x_1$.
Set
\[
b_1=b_1'+y_1 \circ \pi': \bigoplus_{i=0}^N(M^{\otimes_Ri}\otimes_RP) \to M\otimes_RQ.
\]
Then one gets that
\begin{align*}
b_1 \circ \widetilde{a}
&=(b_1'+y_1 \circ \pi') \circ \widetilde{a}       & b_1 \circ c_P &=(b_1'+y_1 \circ \pi') \circ c_P\\
&=b_1' \circ \widetilde{a}+y_1 \circ \pi'
\circ \widetilde{a}                               &&=b_1' \circ c_P+y_1 \circ \pi' \circ c_P \\
&=b_1' \circ \widetilde{a}+y_1 \circ f \circ \pi
\text{  \quad\quad\quad         and }
                                                  &&=M\otimes b_0.\\
&=b_1' \circ \widetilde{a}+x_1 \circ \pi
                                  \\
&=b_1' \circ \widetilde{a}+(l_1-b_1'\circ \widetilde{a})\\
&=l_1
\end{align*}
Iterating the above argument,
one gets an $R$-homomorphism
\[
b_j:\bigoplus_{i=0}^N(M^{\otimes_Ri}\otimes_RP)
\to M^{\otimes_Rj}\otimes_RQ\]
such that $b_j \circ \widetilde{a}=l_j$ and $b_j \circ c_P=M\otimes b_{j-1}$
for $2\leqslant j \leqslant N$. Thus, there exits a $T_R(M)$-homomorphism
\[
\widetilde{b}:=(b_0,b_1, \cdots, b_N)^\tr: \Ind(P) \to \Ind(Q)
\]
such that $\widetilde{b} \circ \widetilde{a}=\widetilde{l}$, and so $\Hom_{T_R(M)}(\widetilde{a},\Ind(Q))$ is epimorphic, as desired.

\textbf{For (2)}:
It is clear that there exists a short exact sequence
\[
\xi: \quad 0 \to (Z,w) \overset{r}\longrightarrow \Ind(P')
             \overset{h}\longrightarrow (X,u) \to 0
\]
of $T_R(M)$-modules with $P' \in \Proj(R)$.
Next, we show that $(Z,w) \in \Phi(\GP(R))$.
By the Snake Lemma,
there exists a short exact sequence
\[0 \to \cok(w) \to P' \to \cok(u) \to 0\]
of $R$-modules. We mention that $\cok(u)\in \GP(R)$. Then one gets that $\cok(w) \in \GP(R)$ as $\GP(R)$ is closed under taking kernels of epimorphisms.
On the other hand, $\Tor^R_1(M, M^{\otimes_Ri}\otimes_R\cok(u))=0$ for $0\leqslant i\leqslant N$ by the assumption, so it follows from Lemma \ref{lem of Tor vanish} that $\Tor^R_1(M,X)=0$. Thus, $M\otimes r$ is a monomorphism, which yields that $w$ is a monomorphism.
Hence, $(Z,w) \in \Phi(\GP(R))$.

Then by an argument similar to the proof of (1),
one gets that the sequence $\Hom_{T_R(M)}(\xi,\Ind(Q))$ is exact.
\end{prf*}

As an immediate consequence of
Lemma \ref{GPTRM in Phi} and Proposition \ref{Phi in GPTRM},
we have the next result.

\begin{thm}\label{THM GP com}
Suppose that the $R$-bimodule $M$ is compatible and admissible.
Then there is an equality
$\GP(T_R(M)) = \Phi(\GP(R))$.
\end{thm}

Now we have the next corollary by Proposition \ref{example of admissible and compatible}.

\begin{cor}\label{G-proj}
Suppose that the $R$-bimodule $M$ satisfies the condition $(\sfT)$, and $\pd_RM<\infty$ and $\fd_{\Rop}M<\infty$. Then there is an equality
$\GP(T_R(M)) = \Phi(\GP(R))$.
\end{cor}

The following example is due to Cui, Rong and Zhang \exacite[4.7]{2022ZHANG},
which shows that there exists a bimodule $M$ over some algebra satisfying
all the conditions in the above corollary.

\begin{exa} \label{example of qviver}
Let $k$ be a field and $kQ$ the path algebra
associated to the quiver
\begin{align*}
\xymatrix{
Q:\  1 \ar[r]_{} & 2 \ar[r]_{} & 3 \ar[r]_{}
  & \cdots \ar[r]_{} & n \ar@/_3ex/[llll]_{}.}
\end{align*}
Suppose that $J$ is the ideal of $kQ$ generated by all the arrows.
Then $R=kQ/J^h$ is a self-injective algebra for $2\leqslant h \leqslant n$.
Denote by $e_i$ the idempotent element corresponding to the vertex $i$.
Then one has $e_jRe_i = 0$ whenever
$1\leqslant i< j\leqslant n$ and $j-i \geqslant h$.
Let $M = Re_i\otimes_k e_jR$. Then $M$ is an $R$-bimodule and projective on both sides, and $M\otimes_RM\cong Re_i\otimes_k (e_jR\otimes_RRe_i)\otimes_k e_jR=0$.
\end{exa}

We mention that the definitions and results in the former part of this section
have dual versions, which are listed in the following.

\begin{dfn} \label{dfn of cocompatible}
An $R$-bimodule $X$ is said to be \emph{co}-\emph{compatible} if
the following two conditions hold
for a complete injective resolution $I^\bullet$ of $T_R(M)^{\op}$-modules:
\begin{prt}
\item[(C1')] The complex $\Hom_{\Rop}(X,I^\bullet)$ is exact.
\item[(C2')] The complex $\Hom_{T_R(M)^{\op}}(\Coind(\Hom_{\Rop}(X,E)),I^\bullet)$ is exact for each $E\in\Inj(\Rop)$.
\end{prt}
\end{dfn}

\begin{dfn} \label{dfn of coadmissible}
An $R$-bimodule $X$ is said to be \emph{co}-\emph{admissible} if
\[
\Ext_{\Rop}^1(\Hom_{\Rop}(X^{\otimes_Ri},E),H)=0=\Ext_{\Rop}^1(X, \Hom_{\Rop}(X^{\otimes_Ri},H))
\]
for each $H \in \GI(\Rop)$, $E \in \Inj(\Rop)$ and $i\geqslant 0$.
\end{dfn}

\begin{lem} \label{dual lemma4.2 for (T)}
Let $X$ be an $R$-bimodule satisfying the condition $(\sfT)$.
Then one has $\Ext^{\geqslant 1}_{\Rop}(X,\Hom_{\Rop}(X^{\otimes_Ri},E))=0$ for each $E\in \Inj(\Rop)$ and $i \geqslant 1$.
In this case, the following statements are equivalent for every $\Rop$-module $Y$.
\begin{eqc}
\item $\Ext_{\Rop}^{\geqslant 1}(X,\Hom_{\Rop}(X^{\otimes_Rj},Y))=0$ for any $j \geqslant 0$.
\item $\Ext_{\Rop}^{\geqslant 1}(X^{\otimes_Rs},\Hom_{\Rop}(X^{\otimes_Rj},Y))=0$
for any $s \geqslant 1$ and $j \geqslant 0$.
\item $\Ext_{\Rop}^{\geqslant 1}(X^{\otimes_Rs},Y)=0$ for any $s \geqslant 1$.
\end{eqc}
\end{lem}
\begin{prf*}
By the assumption, one has $\Tor_{\geqslant 1}^R(X, X^{\otimes_Ri})=0$ for each $i\geqslant 1$, and so $\Ext^{\geqslant 1}_{\Rop}(X,\Hom_{\Rop}(X^{\otimes_Ri},E))\is\Hom_{\Rop}(\Tor_{\geqslant 1}^R(X, X^{\otimes_Ri}),E)=0$ for each $E \in \Inj(\Rop)$.
Then by a dual argument of the proof of Lemma \ref{lemma4.2 for (T)},
one can prove the equivalence of $(i)$-$(iii)$.
\end{prf*}

\begin{prp} \label{example of coadmissible and cocompatible}
Suppose that the $R$-bimodule $M$ satisfies the condition $(\sfT)$.
If $\fd_{R}M<\infty$ and $\pd_{\Rop}M<\infty$,
then $M$ is co-compatible and co-admissible.
\end{prp}
\begin{prf*}
We first prove that $M$ is co-compatible. Let
$$I^\bullet:\ \xymatrix@C=0.5cm{
\cdots \to \Coind(E^{-1}) \ar[r]^{\quad \,\,\,\,\, d^{-1}} & \Coind(E^0)
\ar[r]^{d^0 \quad \,\,} & \Coind(E^1) \to \cdots}$$
be a complete injective resolution of injective $T_R(M)^{\op}$-modules with each $E^j \in \Inj(\Rop)$. Then one gets an exact complex
$$E^\bullet:\ \xymatrix@C=0.5cm{ \cdots \to \bigoplus_{i=0}^N\Hom_{R^{\op}}(M^{\otimes_Ri},E^{-1})
\ar[r]^{\,\, d^{-1}} & \bigoplus_{i=0}^N\Hom_{R^{\op}}(M^{\otimes_Ri},E^{0}) \to \cdots}$$
of $\Rop$-modules.
To prove the exactness of $\Hom_{\Rop}(M,I^\bullet)$,
it suffices to prove that the complex $\Hom_{\Rop}(M,E^\bullet)$ is exact.
Fix a projective resolution $0 \to P^n \to \cdots \to  P^0 \to M \to 0$
of the $\Rop$-module $M$.
By Lemma \ref{dual lemma4.2 for (T)}, one has the equality $\Ext^{\geqslant 1}_{\Rop}(M,\bigoplus_{i=0}^N\Hom_{\Rop}(M^{\otimes_Ri},E^{j}))=0$,
so there is an exact sequence
$$0 \to \Hom_{\Rop}(M,E^\bullet) \to \Hom_{\Rop}(P^0,E^\bullet)
\to \cdots \to \Hom_{\Rop}(P^n,E^\bullet) \to 0$$
of complexes, which yields that the complex $\Hom_{\Rop}(M,E^\bullet)$ is exact.

On the other hand, it follows from Lemma \ref{dual lemma4.2 for (T)} that
$$\Ext_{\Rop}^{\geqslant 1}(M^{\otimes_Ri},\Hom_{\Rop}(M,E))=0,$$
so $\Ext_{\Rop}^{\geqslant 1}(T_R(M),\Hom_{\Rop}(M,E))=0$. Thus,
by Corollary \ref{Inj of Ind-}(b), one has
\[\id_{T_R(M)^{\op}}\Coind(\Hom_{\Rop}(M,E))=\id_{\Rop}(\Hom_{\Rop}(M,E))<\infty,\]
and hence, the complex $\Hom_{T_R(M)^{\op}}(\Coind(\Hom_{\Rop}(M,E)),I^\bullet)$ is exact. This implies that $M$ is co-compatible.

Next, we prove that $M$ is co-admissible.
Indeed, one has $\fd_{R}(M^{\otimes_Ri})<\infty$ for each $0\leqslant i \leqslant N$ by Lemma \ref{lemma4.5 for (T)}, so $\id_{\Rop}(\Hom_{\Rop}(M^{\otimes_Ri},E))<\infty$
for each $E \in \Inj(\Rop)$,
which implies that $\Ext_{\Rop}^{\geqslant 1}(\Hom_{\Rop}(M^{\otimes_Ri},E),H)=0$
for each $H \in \GI(\Rop)$. On the other hand,
by Lemma \ref{lemma4.5 for (T)} again, one has $\pd_{\Rop}(M^{\otimes_Ri})<\infty$,
which ensures that $\Hom_{\Rop}(M^{\otimes_Ri}, U^\bullet)$ is exact
for a complete injective resolution $U^\bullet$ of $H$.
Therefore, $\Ext_{\Rop}^{\geqslant 1}(M^{\otimes_Ri}, H)=0$.
This yields that
$\Ext^{\geqslant 1}_{\Rop}(M,\Hom_{\Rop}(M^{\otimes_Ri},H))=0$
by Lemma \ref{dual lemma4.2 for (T)}.
Thus, $M$ is co-admissible.
\end{prf*}

The next result is a dual version of Theorem \ref{THM GP com}, which can be proved dually.

\begin{thm} \label{THM GITRM = Psi cons}
Suppose that the $R$-bimodule $M$ is co-compatible and co-admissible.
Then there is an equality
$\GI(T_R(M)^{\op}) = \Psi(\GI(\Rop))$.
\end{thm}

The following result is an immediate consequence of
Proposition \ref{example of coadmissible and cocompatible}.

\begin{cor}\label{G-inj}
Suppose that the $R$-bimodule $M$ satisfies the condition $(\sfT)$, and $\fd_{R}M<\infty$ and $\pd_{\Rop}M<\infty$. Then there is an equality
$\GI(T_R(M)^{\op}) = \Psi(\GI(\Rop))$.
\end{cor}

\section{Coherence of tensor rings and Gorenstein flat $T_R(M)$-modules}
\label{Gorenstein injective right T_R(M)-modules}
\noindent
In this section, we investigate the coherence of $T_R(M)$ and
characterize Gorenstein flat $T_R(M)$-modules.

\subsection{The coherence of tensor rings}
\label{The coherence of T_R(M)}
In this subsection, we give a characterization for the coherence of $T_R(M)$.
Before giving the main result Theorem \ref{coherence} of this subsection,
we finish a few preparatory works.


\begin{lem} \label{submodule of TRM}
Let $(X,u)$ be a $T_R(M)$-module.
Then the following statements are equivalent.
\begin{eqc}
\item
$(K,a)$ is a submodule of $(X,u)$ with the inclusion map $\iota$.
\item
$K$ is a submodule of $X$ with the inclusion map $\iota$,
and $\im(u\circ(M\otimes\iota))\subseteq K$.
\end{eqc}
\end{lem}
\begin{prf*}
\proofofimp{i}{ii} If $(K,a)$ is a submodule of $(X,u)$ with the inclusion map $\iota$,
then $K$ is obvious a submodule of $X$ with $\iota$.
Moreover, one has $a=\iota \circ a=u \circ (M\otimes\iota)$, and so $\im(u \circ (M\otimes\iota))=\im a\subseteq K$.

\proofofimp{ii}{i} Suppose that $K$ is a submodule of $X$ with the inclusion map $\iota$ and $\im(u \circ (M\otimes\iota))\subseteq K$. Set $a=u \circ (M\otimes\iota): M\otimes_RK\to K$. Then one gets a $T_R(M)$-homomorphism $\iota: (K,a) \to (X,u)$, so $(K,a)$ is a submodule of $(X,u)$.
\end{prf*}

\begin{lem} \label{f.p.}
If $(X,u)$ is a finitely presented $T_R(M)$-module,
then $\cok(u)$ is a finitely presented $R$-module.
The converse statement holds whenever $u$ is a monomorphism.
\end{lem}
\begin{prf*}
Suppose that $(X,u)$ is a finitely presented $T_R(M)$-module.
Then there is a short exact sequence
$0 \to (Z,w) \to (Y,v)\to (X,u) \to 0$
of $T_R(M)$-modules, where $(Y,v)$ is
finitely generated and free, and $(Z,w)$ is finitely generated.
By Lemma \ref{projectives in F-Rep}, one has $(Y,v)\cong\Ind(P)$
with $P$ a finitely generated and free $R$-module.
By the Snake Lemma, there exists an exact sequence
$$0 \to \kernel(u) \to \cok(w) \to P\to \cok(u) \to 0$$
of $R$-modules, where $\cok(w)$ is finitely generated
by Lemma \ref{projectives in F-Rep}.
This implies that $\cok(u)$ is finitely presented.

Conversely, suppose that $\cok(u)$ is a finitely presented $R$-module and $u$ is a monomorphism. Then $(X,u)$ is finitely generated by Lemma \ref{projectives in F-Rep}.
Take a short exact sequence
$0 \to (Z,w) \to \Ind(P)\to (X,u) \to 0$
of $T_R(M)$-modules with $P$ a finitely generated and free $R$-module.
To complete the proof, it remains to show that $(Z,w)$ is finitely generated.
Indeed, by the Snake Lemma, there exists a short exact sequence
$$0 \to \cok(w) \to P \to \cok(u) \to 0$$
of $R$-modules.
Since $\cok(u)$ is finitely presented, one gets that $\cok(w)$ is finitely generated,
and so $(Z,w)$ is finitely generated by Lemma \ref{projectives in F-Rep} again.
\end{prf*}

\begin{lem} \label{coherent TRM module}
Consider the following statements:
\begin{eqc}
\item $T_R(M)$ is left coherent.
\item For each submodule $(K, a)$ of $\Ind(R)$, if $K/\im a$ is a finitely generated $R$-module, then it is finitely presented.
\end{eqc}
Then \eqclbl{i}$\!\implies\!$\eqclbl{ii} holds true. Furthermore, if $M$ is flat as an $\Rop$-module, then the two statements are equivalent.
\end{lem}
\begin{prf*}
Note that the tensor ring $T_R(M)$, as a $T_R(M)$-module, coincides with
$\Ind(R)=(\bigoplus_{i=0}^NM^{\otimes_Ri},c_R)$, where $c_R$ is a monomorphism.
This implies that $T_R(M)$ is left coherent if and only if
every finitely generated submodule $(K, a)$ of $\Ind(R)$ is finitely presented.

\proofofimp{i}{ii} Suppose that $K/\im a$ is a finitely generated $R$-module. Then $(K, a)$ is finitely generated by Lemma \ref{projectives in F-Rep}, and hence, $(K, a)$ is finitely presented.
This yields that $K/\im a$ is finitely presented by Lemma \ref{f.p.}.

\proofofimp{ii}{i} Suppose that $(K, a)$ is a finitely generated submodule of $\Ind(R)$ with the inclusion map $\iota$.
Then $K/\im a$ is a finitely generated $R$-module by Lemma \ref{projectives in F-Rep},
and hence, $K/\im a$ is finitely presented.
By Lemma \ref{submodule of TRM}, we see that $a=c_R\circ(M\otimes\iota)$.
Since $M$ is a flat $\Rop$-module by assumption,
it follows that $a$ is a monomorphism.
Thus $(K, a)$ is finitely presented by Lemma \ref{f.p.}. Therefore, $T_R(M)$ is left coherent.
\end{prf*}

\begin{lem} \label{coherent2 TRM module}
The following statements are equivalent.
\begin{eqc}
\item $R$ is left coherent and $M$ is finitely presented as an $R$-module.
\item $M$ is finitely generated as an $R$-module, and for each submodule $(K,a)$ of\; $\Ind(R)$, if $K/\im a$ is a finitely generated $R$-module, then it is finitely presented.
\end{eqc}
\end{lem}
\begin{prf*}
\proofofimp{i}{ii} Let $(K,a)$ be a submodule of $\Ind(R)$ with the inclusion map $\iota$
such that $K/\im a$ is a finitely generated $R$-module.
We then prove that $K/\im a$ is finitely presented.
By Lemma \ref{submodule of TRM}, one gets that $K$ is a submodule of $\bigoplus_{i=0}^NM^{\otimes_Ri}$ with the inclusion map $\iota$,
and $\im(c_R \circ (M\otimes\iota))\subseteq K$.
Then for each $0\leqslant i \leqslant N$,
there exists a submodule $K_i$ of $M^{\otimes_Ri}$ with the inclusion map $\iota_i$,
such that
$$K=K_0\oplus K_1\oplus\cdots \oplus K_N,$$
and
$$\iota=\begin{pmatrix}
\iota_0 & 0 & 0 & \cdots & 0\\
0 & \iota_1 & 0 & \cdots & 0\\
0 & 0 & \iota_2 & \cdots & 0\\
\cdots & \cdots & \cdots & \cdots & \cdots\\
0 & 0 & 0 & \cdots & \iota_N\\
\end{pmatrix}$$
$$a=\begin{pmatrix}
0 & 0 & 0 & \cdots & 0 & 0\\
M\otimes\iota_0 & 0 & 0 & \cdots & 0 & 0\\
0 & M\otimes\iota_1 & 0 & \cdots & 0 & 0\\
\cdots & \cdots & \cdots & \cdots & \cdots & \cdots\\
0 & 0 & 0 & \cdots & M\otimes\iota_{N-1} & 0\\
\end{pmatrix}.
$$
Moreover, one has $K/\im a \is K_0\oplus K_1/\im(M\otimes\iota_0) \oplus \cdots \oplus K_N/\im(M\otimes\iota_{N-1})$.
Since $K/\im a$ is finitely generated, one gets that $K_0$ and $K_j/\im(M\otimes\iota_{j-1})$ are finitely generated $R$-modules for each $1\leqslant j\leqslant N$.
It is clear that $K_0$ is finitely presented as $R$ is left coherent.
We then prove that $K_j/\im(M\otimes\iota_{j-1})$ is finitely presented for each $1\leqslant j\leqslant N$.

Consider the short exact sequence
$$0 \to \im(M\otimes\iota_{0}) \to K_1 \to K_1/\im(M\otimes\iota_{0}) \to 0$$
of $R$-modules.
The $R$-module $M\otimes_RK_0$ is finitely generated,
so $\im(M\otimes\iota_{0})$ is finitely generated,
and hence, $K_1$ is a finitely generated submodule of $M$.
This yields that $K_1$ is finitely presented as $R$ is left coherent
and $M$ is a finitely presented $R$-module,
and so $K_1/\im(M\otimes\iota_{0})$ is finitely presented; see \cite[25.1]{1991WRobert}.
Consider the short exact sequence
$$0 \to \im(M\otimes\iota_{1}) \to K_2 \to K_2/\im(M\otimes\iota_{1}) \to 0.$$
Then $\im(M\otimes\iota_{1})$ is finitely generated as so is $M\otimes_RK_1$,
and so $K_2$ is a finitely generated submodule of $M\otimes_RM$.
Note that $M$ is a finitely presented $R$-module.
For a family $\{X_i\}_{i \in I}$ of $R$-modules, one can easily to check that
$\prod_{i \in I} X_i \otimes_R(M\otimes_RM)\cong (\prod_{i \in I} X_i \otimes_RM)\otimes_RM
\cong \prod_{i \in I}(X_i \otimes_RM)\otimes_RM \cong \prod_{i \in I}((X_i \otimes_RM)\otimes_RM)
\cong \prod_{i \in I}(X_i \otimes_R(M\otimes_RM))$.
This implies that $M\otimes_RM$ is a finitely presented $R$-module.
Then $K_2$ is finitely presented as $R$ is left coherent,
and hence, $K_2/\im(M\otimes\iota_{1})$ is finitely presented by \cite[25.1]{1991WRobert} again.
Iterating the above arguments, one gets that the $R$-modules $K_j/\im(M\otimes\iota_{j-1})$ is finitely presented for each $1\leqslant j\leqslant N$, as desired.
Therefore, $K/\im a$ is finitely presented.

\proofofimp{ii}{i} We first prove that $R$ is left coherent. To this end,
let $K_0$ be a finitely generated left ideal of $R$ with the inclusion map $\iota_0$,
then it suffices to show that $K_0$ is finitely presented. We mention that
$(K_0\oplus M\oplus\cdots \oplus M^{\otimes_RN},b)$ is a submodule of $\Ind(R)$, where
\[b=\begin{pmatrix}
0 & 0 & 0 & \cdots & 0\\
M\otimes\iota_{0} & 0 & 0 & \cdots & 0\\
0 & 1 & 0 & \cdots & 0\\
\cdots & \cdots & \cdots & \cdots & \cdots\\
0 & 0 & 0 & \cdots & 1\\
\end{pmatrix}_{(N+1)\times N}. \]
It is easy to check that
$(K_0\oplus M\oplus\cdots \oplus M^{\otimes_RN})/\im b\is K_0\oplus M/\im(M\otimes\iota_{0})$.
Since $M$ is a finitely generated $R$-module, one gets that $(K_0\oplus M\oplus\cdots \oplus M^{\otimes_RN})/\im b$ is finitely generated, and so it is finitely presented by the assumption. Thus, $K_0$ is finitely presented, as desired.

Next, we prove that $M$ is a finitely presented $R$-module.
We mention that
$(M\oplus M^{\otimes_R2}\oplus\cdots \oplus M^{\otimes_RN},d)$ is a submodule of $\Ind(R)$, where
\[d=\begin{pmatrix}
0 & 0 & 0 & \cdots & 0\\
1 & 0 & 0 & \cdots & 0\\
0 & 1 & 0 & \cdots & 0\\
\cdots & \cdots & \cdots & \cdots & \cdots\\
0 & 0 & 0 & \cdots & 1\\
\end{pmatrix}_{(N,N-1)}. \]
It is easy to check that $(M\oplus M^{\otimes_R2}\oplus\cdots \oplus M^{\otimes_RN})/\im d \is M$. Since $M$ is a finitely generated $R$-module, one gets that $(M\oplus M^{\otimes_R2}\oplus\cdots \oplus M^{\otimes_RN})/\im d $ is finitely generated, and so it is finitely presented by the assumption. Thus, $M$ is finitely presented.
\end{prf*}

As an immediate consequence of Lemmas \ref{coherent TRM module} and
\ref{coherent2 TRM module}, we obtain the following result, which gives a characterization for the coherence of $T_R(M)$.

\begin{thm}\label{coherence}
Consider the following statements:
\begin{eqc}
\item $M$ is finitely generated as an $R$-module and $T_R(M)$ is left coherent.
\item $M$ is finitely presented as an $R$-module and $R$ is left coherent.
\end{eqc}
Then \eqclbl{i}$\!\implies\!$\eqclbl{ii} holds true. Furthermore, if $M$ is flat as an $\Rop$-module, then the two statements are equivalent.
\end{thm}

\subsection{Gorenstein flat $T_R(M)$-modules}\label{Gorenstein flat T_R(M)-modules}
In this subsection, we characterize Gorenstein flat $T_R(M)$-modules, and give the proof of Theorem \ref{THM GP}(3) in the introduction.

Recall from \dfncite[3.1]{GHD} that an exact sequence
$$\xymatrix@C=0.5cm{
F^\bullet: \cdots \to F^{-1} \ar[r]^{\qquad \quad d^{-1} \,} & F^0
\ar[r]^{\,\, d^0 \qquad } & F^1 \to \cdots,}$$
of flat $R$-modules is called \emph{complete flat resolution} if it remains exact after applying the functor $E\otimes_R-$ for every injective $\Rop$-module $E$.
An $R$-module $X$ is called \emph{Gorenstein flat}
provided that there exists a complete flat resolution $F^\bullet$ such that
$X \cong \kernel(F^0 \to F^1)$.
We denote by $\GF(R)$ the subcategory of Gorenstein flat $R$-modules.

According to \thmcite[3.6]{GHD}, if $G\in\GF(R)$ then $G^+\in\GI(\Rop)$, and the converse holds true whenever $R$ is right coherent. Relying on this fact, we have the following result, which can be proved easily.

\begin{lem}\label{Phi GF=Psi GI}
Let $(X,u)$ be a $T_R(M)$-module.
If $(X,u)\in\Phi(\GF(R))$,
then $(X,u)^+\in\Psi(\GI(\Rop))$.
Furthermore, the converse holds true whenever $R$ is right coherent.
\end{lem}

\begin{lem}\label{Gf-contain}
Let $R$ be a right coherent ring.
Suppose that the $R$-bimodule $M$ satisfies the condition $(\sfT)$,
and $\fd_{R}M<\infty$ and $\pd_{\Rop}M<\infty$.
Then there is an inclusion $\GF(T_R(M))\subseteq\Phi(\GF(R))$.
\end{lem}
\begin{prf*}
Let $(X,u)$ be a Gorenstein flat $T_R(M)$-module. Then $(X,u)^+$ is Gorenstein injective, and so it is in $\Psi(\GI(\Rop))$ by Corollary \ref{G-inj}. Thus, one has $(X,u)\in\Phi(\GF(R))$ by Lemma \ref{Phi GF=Psi GI}.
\end{prf*}

The next result gives a characterization for Gorensein flat $T_R(M)$-modules.

\begin{thm}\label{G-flat}
Let $R$ be a right coherent ring.
Suppose that $M$ is flat as an $R$-module and finitely presented as an $\Rop$-module, and $\pd_{\Rop}M<\infty$.
Then there is an equality $\GF(T_R(M)) = \Phi(\GF(R))$.
\end{thm}
\begin{prf*}
Since $M$ is flat as an $R$-module, one gets that $M^{\otimes_Ri}\otimes_RP$ is a flat  $R$-module for each $P\in\Proj(R)$ and every $i\geqslant 1$.
Hence, the condition $(\sfT)$ holds clearly.
Thus, by Lemma \ref{Gf-contain}, one has $\GF(T_R(M))\subseteq\Phi(\GF(R))$. In the following, we prove the inclusion $\Phi(\GF(R))\subseteq\GF(T_R(M))$. To this end, let $(X,u)$ be in $\Phi(\GF(R))$. Then one has $(X,u)^+\in\GI(T_R(M)^{\op})$ by Corollary \ref{G-inj} and Lemma \ref{Phi GF=Psi GI}, and so $(X,u)\in\GF(T_R(M))$ as $T_R(M)$ is right coherent by Theorem \ref{coherence}.
\end{prf*}

\begin{rmk}
So far, we do not know whether the condition that $R$ is right coherent in Theorem \ref{G-flat} can be removed. Recall from \cite{2020Jan} that an $R$-module $X$ is called
\emph{projectively coresolved Gorenstein flat} if there is an exact sequence
$$\xymatrix@C=0.5cm{
P^\bullet: \cdots \to P^{-1} \ar[r]^{\qquad \quad d^{-1} \,} & P^0
\ar[r]^{\,\, d^0 \qquad } & P^1 \to \cdots,}$$
of projective $R$-modules with $X \cong \kernel(P^0 \to P^1)$,
such that it remains exact after applying the functor
$E\otimes_R-$ for every injective $\Rop$-module $E$. We do not know whether Theorem \ref{G-flat} holds for projectively coresolved Gorenstein flat modules.
\end{rmk}

\section{Applications}
\label{Applications}
\noindent
In this section, we give some applications to
trivial ring extensions and Morita context rings.

\subsection{The trivial extension of rings}
\label{The trivial extension of rings}
\noindent
Let $R$ be a ring and $M$ an $R$-bimodule.
There exists a ring $R\ltimes M$,
where the addition is componentwise
and the multiplication is given by
$(r_1, m_1)(r_2, m_2) = (r_1r_2,r_1m_2 + m_1r_2)$ for $r_1, r_2 \in R$ and $m_1, m_2 \in M$.
This ring is called the \emph{trivial extension} of
the ring $R$ by the $R$-bimodule $M$; see \cite{TRIEXT1975} and \cite{TRIEXT1971}.

Suppose that the $R$-bimodule $M$ is $1$-nilpotent,
that is, $M\otimes_RM=0$.
Then it is easy to see that the tensor ring $T_R(M)$
is nothing but the trivial ring extension $R\ltimes M$.
One can immediately get the following results by
Theorem \ref{THM GP}.

\begin{cor} \label{GPTRM = Phi in trivial extension}
Suppose that $M$ is a $1$-nilpotent $R$-bimodule.
\begin{rqm}
\item
If $\Tor_{\geqslant 1}^R(M,M\otimes_RP)=0$ for each $P\in \Proj(R)$,
$\pd_{R}M<\infty$ and $\fd_{\Rop}M<\infty$, then there is an equality $\GP(R\ltimes M) = \Phi(\GP(R))$.

\item
If $\Tor_{\geqslant 1}^R(M,M\otimes_RP)=0$ for each $P\in \Proj(R)$, $\fd_{R}M<\infty$ and $\pd_{\Rop}M<\infty$, then there is an equality $\GI((R\ltimes M)^{\op}) = \Psi(\GI(\Rop))$.

\item
If $R$ is  right coherent, $M$ is flat as an $R$-module and finitely presented as an $\Rop$-module, and $\pd_{\Rop}M<\infty$, then there is an equality $\GF(R\ltimes M)=\Phi(\GF(R))$.
\end{rqm}
\end{cor}

\begin{rmk} \label{GP in trivial extension}
Suppose that $M$ is a $1$-nilpotent $R$-bimodule and the conditions in Corollary \ref{GPTRM = Phi in trivial extension}(1) hold. Recall from \ref{The stalk functor and its adjoints} the definition of the stalk functor $S: \Mod(R)\to \Mod(R\ltimes M)$. Fix a projective resolution
$$0 \to P^n \to \cdots \to P^1 \to P^0 \to M \to 0$$
of the $R$-module $M$. Since $\Tor_{\geqslant 1}^R(M,M)=0$ and $M\otimes_RM=0$ by the assumption, there is an exact sequence
$$0 \to M\otimes_RP^n \to \cdots \to M\otimes_RP^1 \to M\otimes_RP^0 \to 0.$$
Thus, one gets a projective resolution of $S(M)$ as follows:
$$0 \to (P^n\oplus (M\otimes_RP^n), (0,1)^{\tr}) \to
\cdots \to (P^0\oplus (M\otimes_RP^0), (0,1)^{\tr})
\to S(M) \to 0,$$
which yields that $\pd_{R\ltimes M}S(M)<\infty$. We mention that there is an exact sequence
$$0 \to S(M) \to (R\oplus M, (0,1)^{\tr}) \to S(R) \to 0$$
of $R\ltimes M$-modules. Then one has $\pd_{R\ltimes M}S(R)<\infty$ as $(R\oplus M, (0,1)^{\tr})$ is projective. Similarly, one gets that $\fd_{(R\ltimes M)^{\op}}S(R)<\infty$. Thus, $S(R)$ is a generalized compatible $R\ltimes M$-bimodule in the sense of Mao; see \exacite[3.2(1)]{TriExtGPMao}. In this case, Corollary \ref{GPTRM = Phi in trivial extension}(1) is a special case of \corcite[3.6(1)]{TriExtGPMao}. Similarly, the statements (2) and (3) in Corollary \ref{GPTRM = Phi in trivial extension} are special cases of \corcite[4.6(1)]{TriExtGPMao} and \corcite[4.8(1)]{TriExtGPMao}, respectively.

We mention that the above argument actually shows that if $M$ is a $1$-nilpotent $R$-bimodule such that $\Tor_{\geqslant 1}^R(M,M)=0$, $\pd_{R}M<\infty$ and $\fd_{\Rop}M<\infty$, then $S(R)$ is a generalized compatible $R\ltimes M$-bimodule. Using this fact, one can remove the assumption that $R$ is a perfect ring in \exacite[3.7]{TriExtGPMao}.
\end{rmk}

Mao studied the coherence of a trivial ring extension $R\ltimes M$;
see \prpcite[2.3]{2023MaoCoherence}.
In the following, we give a new characterization,
which is an immediate consequence of Theorem \ref{coherence}.

\begin{cor} \label{coherence in trivial extension}
Suppose that $M$ is a $1$-nilpotent $R$-bimodule. Consider the following statements:
\begin{eqc}
\item $M$ is finitely generated as an $R$-module and $R\ltimes M$ is left coherent.
\item $M$ is finitely presented as an $R$-module and $R$ is left coherent.
\end{eqc}
Then \eqclbl{i}$\!\implies\!$\eqclbl{ii} holds true. Furthermore, if $M$ is flat as an $\Rop$-module, then the two statements are equivalent.
\end{cor}

\subsection{Morita context rings}
\label{Morita context rings}
Let $A$ and $B$ be two rings,
and let $_AV_B$ and $_BU_A$ be two bimodules,
$\phi : U\otimes_AV \to B$ a homomorphism of $B$-bimodules,
and $\psi : V\otimes_BU \to A$ a homomorphism of $A$-bimodules.
We assume further that $\phi(u\otimes v)u'=u\psi(v \otimes u')$
and $v\phi(u\otimes v')=\psi(v \otimes u)v'$
for all $u, u' \in U$ and $v, v' \in V$.
Associated with a \emph{Morita context} $(A,B,U,V,\phi,\psi)$,
there exists a \emph{Morita context ring}
\[\Lambda_{(\phi,\psi)}=\begin{pmatrix}A & V \\ U & B\end{pmatrix},\]
where the addition of elements is componentwise and multiplication is given by
\[\begin{pmatrix}a & v \\ u & b\end{pmatrix}\cdot
\begin{pmatrix}a' & v' \\ u' & b'\end{pmatrix}=
\begin{pmatrix}aa'+\psi(v \otimes u') & av'+vb' \\
ua'+bu' & bb'+ \phi(u \otimes v') \end{pmatrix};\]
see \cite{BASS1968, MORITA1958} for more details.
Following \thmcite[1.5]{GF1982}, one can view a $\Lambda_{(\phi,\psi)}$-module
as a quadruple $(X,Y,f,g)$ with $X \in \Mod(A)$, $Y \in \Mod(B)$,
$f \in \Hom_B(U\otimes_AX,Y)$, and $g \in \Hom_A(V\otimes_BY,X)$. Also, one can view a $\Lambda_{(\phi,\psi)}^{\op}$-module as a quadruple
$[W,N,s,t]$ with $W \in \Mod(A^{\op})$, $N \in \Mod(B^{\op})$,
$s \in \Hom_{B^{\op}}(N,\Hom_{A^{\op}}(U,W))$, and $t \in \Hom_{A^{\op}}(W,\Hom_{B^{\op}}(V,N))$.

It follows from \prpcite[2.5]{GFARTIN} that Morita context rings are trivial ring extensions
whenever both $\phi$ and $\psi$ are zero.
More precisely, consider the Morita context ring $\Lambda_{(0,0)}$.
There exists an isomorphism of rings:
$$\Lambda_{(0,0)}\overset{\cong}\longrightarrow (A\times B)\ltimes (U \oplus V)
\,\,\text{via}\,\, \begin{pmatrix}a & v \\ u & b\end{pmatrix} \mapsto
((a,b),(u,v)).$$
Thus, there exists an isomorphic functor
$$\mu: \Mod(\Lambda_{(0,0)}) \to \Mod((A\times B)\ltimes (U \oplus V))\ \mathrm{via}\
(X,Y,f,g) \mapsto ((X,Y),(g,f)),$$
where $(g,f)$ is from
$(U \oplus V)\otimes_{A\times B}(X, Y) \cong (V\otimes_BY,U\otimes_AX)$ to $(X, Y)$.

We mention that $(U \oplus V)\otimes_{A\times B}(U \oplus V) \cong (U\otimes_AV) \oplus (V\otimes_BU)$.
Then the $A\times B$-bimodule $U \oplus V$ is $1$-nilpotent if and only if $U\otimes_AV=0=V\otimes_BU$.
Thus, we obtain the next result by Corollary \ref{GPTRM = Phi in trivial extension}.

\begin{cor} \label{con GP in morita ring}
Let $\Lambda_{(0,0)}$ be a Morita context ring with $U\otimes_AV=0=V\otimes_BU$, and let $(X,Y,f,g)$ be a $\Lambda_{(0,0)}$-module and $[W,N,s,t]$ a $\Lambda_{(0,0)}^{\op}$-module.
\begin{rqm}
\item
Suppose that $\Tor^B_{\geqslant 1}(V, U\otimes_AP_1)=0=\Tor^A_{\geqslant 1}(U, V\otimes_BP_2)$ for each $P_1 \in \Proj(A)$ and $P_2 \in \Proj(B)$. If $\pd_{B}U<\infty$, $\fd_{A^{\op}}U<\infty$, $\pd_{A}V<\infty$ and $\fd_{B^{\op}}V< \infty$, then $(X,Y,f,g)\in\GP(\Lambda_{(0,0)})$ if and only if
both $f$ and $g$ are monomorphisms, and $\cok(f)\in\GP(B)$ and $\cok(g)\in\GP(A)$.
\item
Suppose that $\Tor^B_{\geqslant 1}(V, U\otimes_AP_1)=0=\Tor^A_{\geqslant 1}(U, V\otimes_BP_2)$ for each $P_1 \in \Proj(A)$ and $P_2 \in \Proj(B)$. If $\fd_{B}U<\infty$, $\pd_{A^{\op}}U<\infty$, $\fd_{A}V<\infty$ and $\pd_{B^{\op}}V<\infty$, then $[W,N,s,t]\in\GI(\Lambda_{(0,0)}^{\op})$ if and only if both $s$ and $t$ are epimorphisms, and $\kernel(s)\in\GI(B^{\op})$ and $\kernel(t)\in\GI(A^{\op})$.
\item
Suppose that $A$ and $B$ are right coherent rings, $U$ is flat as a $B$-module and finitely presented as an $A^{\op}$-module, and $V$ is flat as an $A$-module and finitely presented as a $B^{\op}$-module. If $\pd_{A^{\op}}U<\infty$ and  $\pd_{B^{\op}}V<\infty$, then $(X,Y,f,g)\in\GF(\Lambda_{(0,0)})$ if and only if both $f$ and $g$ are monomorphisms, and $\cok(f)\in\GF(B)$ and $\cok(g)\in\GF(A)$.
\end{rqm}
\end{cor}
\begin{prf*}
We only prove (1); the statements (2) and (3) can be proved similarly.

By the assumption one can check that
$\Tor^{A\times B}_{\geqslant 1}
(U \oplus V, (U \oplus V)\otimes_{A\times B}(P_1,P_2))=0$ for each
$(P_1,P_2) \in \Proj(A\times B)$, and both $\pd_{A\times B}(U \oplus V)$ and $\fd_{{(A\times B)}^{\op}}(U \oplus V)$ are finite. Thus, it follows from Corollary \ref{GPTRM = Phi in trivial extension}(1) that there is an equality
$$
\GP(\Lambda_{(0,0)})
\cong \GP((A\times B)\ltimes (U \oplus V))
=\Phi(\GP(A\times B)),$$
where an object $((X,Y),(g,f))\in\Mod((A\times B)\ltimes (U \oplus V))$ is in $\Phi(\GP(A\times B))$ if and only if $(g,f)$ is a monomorphism and $\cok(g,f)$ is in $\GP(A\times B)$. This yields that $(X,Y,f,g)\in\GP(\Lambda_{(0,0)})$ if and only if
both $f$ and $g$ are monomorphisms, and $\cok(f)\in\GP(B)$ and $\cok(g)\in\GP(A)$.
\end{prf*}

\begin{rmk}
By Corollary \ref{con GP in morita ring}(1), one gets that if $\pd_{B}U<\infty$, $\pd_{A}V<\infty$, and $U$ and $V$ are flat as an $A^{\op}$-module and as a $B^{\op}$-module, respectively, then $(X,Y,f,g)\in\GP(\Lambda_{(0,0)})$ if and only if
both $f$ and $g$ are monomorphisms and $\cok f\in\GP(B)$ and $\cok g\in\GP(A)$. This fact improves \thmcite[4.2]{LMY25} by removing the condition that $V\otimes_B\GP(B)\subseteq\GP(A)$, $U\otimes_A\GP(A)\subseteq\GP(B)$, $V\otimes_B\Proj(B)\subseteq\Proj(A)$ and $U\otimes_A\Proj(A)\subseteq\Proj(B)$.
\end{rmk}

\begin{exa}\label{exa of motita UV=0}
Suppose that the ring $R$ and the $R$-bimodule $M$ as in
Example \ref{example of qviver}.
Then the Morita context ring
$$\Lambda_{(0,0)}=\begin{pmatrix}R & M \\ M & R\end{pmatrix}$$
satisfies all conditions in Corollary \ref{con GP in morita ring}.
\end{exa}

Yan and Yao characterized the coherence of Morita context rings $\Lambda_{(0,0)}$; see \thmcite[5.4]{GFFP}. In the following, we give a new characterization by Corollary \ref{coherence in trivial extension}.

\begin{cor}\label{coherence in morita context}
Let $\Lambda_{(0,0)}$ be a Morita context ring with $U\otimes_AV=0=V\otimes_BU$.
Consider the following statements:
\begin{eqc}
\item $\Lambda_{(0,0)}$ is left coherent, and $U$ and $V$ are finitely generated as an $A$-module and a $B$-module, respectively.
\item $A$ and $B$ are left coherent, and $U$ and $V$ are finitely presented as an $A$-module and a $B$-module, respectively.
\end{eqc}
Then \eqclbl{i}$\!\implies\!$\eqclbl{ii} holds true. Furthermore, if $U$ and $V$ are flat as a $B^{\op}$-module and an $A^{\op}$-module respectively, then the two statements are equivalent.
\end{cor}


\section*{Acknowledgments}
\noindent
We thank Rongmin Zhu for helpful discussions related to this work.
We also thank the anonymous referee for pertinent suggestions that improved the exposition.

\bibliographystyle{amsplain-nodash}

\def\cprime{$'$}
  \providecommand{\arxiv}[2][AC]{\mbox{\href{http://arxiv.org/abs/#2}{\sf
  arXiv:#2 [math.#1]}}}
  \providecommand{\oldarxiv}[2][AC]{\mbox{\href{http://arxiv.org/abs/math/#2}{\sf
  arXiv:math/#2
  [math.#1]}}}\providecommand{\MR}[1]{\mbox{\href{http://www.ams.org/mathscinet-getitem?mr=#1}{#1}}}
  \renewcommand{\MR}[1]{\mbox{\href{http://www.ams.org/mathscinet-getitem?mr=#1}{#1}}}
\providecommand{\bysame}{\leavevmode\hbox to3em{\hrulefill}\thinspace}
\providecommand{\MR}{\relax\ifhmode\unskip\space\fi MR }
\providecommand{\MRhref}[2]{%
  \href{http://www.ams.org/mathscinet-getitem?mr=#1}{#2}
}
\providecommand{\href}[2]{#2}

\end{document}